\documentclass{amsart}
\usepackage{amsmath,amssymb,amsthm}
\usepackage{graphicx}

\theoremstyle{plain}
\newtheorem{thm}{Theorem}[section] 

\newtheorem{prop}[thm]{Proposition}

\newtheorem{lem}[thm]{Lemma}
\theoremstyle{definition} 
\newtheorem{defn}[thm]{Definition}

\theoremstyle{remark}
\newtheorem{rem}[thm]{Remark}

\title[Tautness of $F$-regular and $F$-pure surface singularities]{On Tautness of Two-dimensional $F$-regular and $F$-pure rational singularities}
\author{Yuki Tanaka}
\address{Graduate School of Mathematical Science, The University of Tokyo, 3-8-1 Komaba Meguro-ku Tokyo 153-8914, Japan.}
\email{yktanaka@ms.u-tokyo.ac.jp}
\date{}

\begin{document}
\maketitle

\begin{abstract}
The weighted dual graph of a two-dimensional normal singularity $(X, x)$ represents the topological nature of the exceptional locus of its minimal log resolution.
$(X, x)$ and its graph are said to be taut if the singularity can be uniquely determined by the graph.
Laufer \cite{Lau2} gave a complete list of taut singularities over $\mathbb{C}$.
In positive characteristics, taut graphs over $\mathbb{C}$ are not necessarily taut and tautness have been studied only for special cases.
In this paper, we prove the tautness of $F$-regular singularities.
We also discuss the tautness of $F$-pure rational singularities.
\end{abstract}

\section{Introduction}

Throughout this paper, we fix an algebraically closed field $k$ of arbitrary characteristic.
Let $(X, x)$ be a two-dimensional normal singularity over $k$,
	that is, a pair consisting of the spectrum of an two-dimensional normal local ring essentially of finite type over $k$ and its unique closed point.
We say that two such singularities are isomorphic if the completions of their local rings are isomorphic to each other.
A morphism $\pi : Y \rightarrow X$ is called a log resolution if it is a proper birational morphism from an nonsingular surface $Y$ and
	its exceptional locus $\pi ^{- 1} (x)$ is a simple normal crossing divisor.
There exists such $\pi$ isomorphic over $X \setminus \{x\}$. (See \cite{Lip}).
Contracting all $(- 1)$-curves with $2$ or less intersections with other components, we obtain a unique minimal log resolution of $(X, x)$.
For the minimal log resolution $\pi : Y \rightarrow X$,
	let $E := \pi^{-1}(x) = \cup_{i=1}^n E_i\subset Y$ be the exceptional locus of $\pi$.
By assumption, each component $E_i$ is a nonsingular projective curve embedded into a nonsingular surface.

\begin{defn}
For an exceptional divisor $Z = \sum_{i=1}^n \nu_i E_i$ supported on $E$, define the \textit{weighted dual graph $\Gamma_Z$ associated to the divisor $Z$} as follows:
\begin{enumerate}
	\item Each irreducible component $E_i$ corresponds to a vertex $v_i$.
	\item An intersection of $E_i$ and $E_j$ corresponds to an edge between $v_i$ and  $v_j$.
		Consequently, there are $E_i \cdot E_j = \#\{ E_i \cap E_j\}$ edges between $v_i$ and  $v_j$.
	\item Each vertex $v_i$ is associated with three integers, the arithmetic genus $p_a(E_i)$, the self-intersection number $- b_i = E_i^2$ and the multiplicity $\nu_i$.
\end{enumerate}
We define the \textit{weighted dual graph associated to the singularity $(X, x)$} by $\Gamma_{X, x} := \Gamma_E$.
Two weighted dual graphs are said to be \textit{isomorphic} to each other if there exists an isomorphism of graphs preserving all corresponding weights simultaneously.
\end{defn}

Note that $\Gamma_{X, x}$ is isomorphic to $\Gamma_{X', x'}$ if $(X, x)$ is isomorphic to $(X', x')$.
Now we give the definition of tautness for two-dimensional normal singularities.

\begin{defn}
Let $(X, x)$ be a two-dimensional normal singularity over an algebraically closed field $k$.
Then $(X, x)$ is said to be \textit{taut} if the following condition is satisfied:
	if $(X', x')$ is another two-dimensional normal singularity over $k$ and $\Gamma_{X', x'}$ is isomorphic to $\Gamma_{X, x}$, then $(X' , x')$ is isomorphic to $(X, x)$.
We say $\Gamma_{X, x}$ is \textit{taut} if $(X, x)$ is taut.
\end{defn}

Laufer \cite{Lau2} gave a complete list of taut singularities over the complex number field $\mathbb{C}$ using the deformation theory of analytic spaces.
In positive characteristics, the classification of taut singularities are far from complete.
Sh\"{u}ller \cite{Sch} recently proved that modulo $p$ reduction of a two-dimensional taut singularity $(X, x)$ over $\mathbb{C}$ is taut for sufficiently large $p$.
In his proof, he did not give a sharp estimation of the characteristics in which the tautness holds.
Lee and Nakayama \cite{LN} proved in arbitrary characteristics that $\Gamma_{X, x}$ is taut if it is a chain with all genera zero.
Artin's list of rational double points (RDP) \cite{Art} tells us that there are both taut RDPs and non-taut RDPs in positive characteristics

$F$-singularities are important classes of singularities in positive characteristics.
As one of these singularities, $F$-regular singularities were introduced by Hochster and Huneke \cite{HH} in the theory of tight closure.
They can be regarded as a characteristic $p > 0$ analogue of log terminal singularities
	because log terminal singularities over $\mathbb{C}$ become $F$-regular after reduction to modulo $p >> 0$ \cite{Har2} \cite{MS} \cite{Smi}.
On the other hand, observing Laufer's list, we can see log terminal singularities over $\mathbb{C}$ are all taut.
Therefore it is natural to ask whether $F$-regular singularities are taut or not.
Our main result gives an affirmative answer to this question.

\begin{thm}
	Every two-dimensional $F$-regular singularity over an algebraically closed field of positive characteristic is taut.
\end{thm}

There is a larger class of $F$-singularity called $F$-pure singularity.
Although $F$-purity is neither a sufficient condition nor a necessary condition to be taut even for a rational singularity,
	there is a relationship between $F$-purity and a kind of ``tautness'' of rational double points.
This is discussed in Section 5.

\section{$F$-singularity and its classification}

We recall the definition of $F$-regular and $F$-pure singularities.

\begin{defn}[\cite{HH}, \cite{HR}]
	Let $(X, x) = (\mathrm{Spec} A, \{\mathfrak{m}\})$ be a two-dimensional normal singularity over an algebraically closed field $k$ of positive characteristic $p$
		and $F: A \rightarrow A$ be the Frobenius endomorphism sending $f \in A$ to $f^p \in A$.
	For each integer $e > 0$, $e$ times iteration of $F$ gives $A$ another $A$-module structure defined by $a \cdot b = a^q b\, (q = p^e)$
		and we denote this module by $F_*^e A$.
	$X$ is said to be \textit{$F$-finite} if $F_* A$ is a finite $A$-module.
	
	Suppose $X$ is $F$-finite.
	\begin{enumerate}
		\item $(X, x)$ is said to be \textit{$F$-regular} if for every $0 \neq c \in A$, there exists an integer $e > 0$
			such that $c F^e : A \rightarrow F_*^e A$ sending $x$ to $c x^{p^e}$ splits as an $A$-module homomorphism.
		\item $(X, x)$ is said to be \textit{$F$-pure} if $F : A \rightarrow F_* A$ splits as an $A$-module homomorphism.
	\end{enumerate}
\end{defn}

$F$-regularity implies $F$-purity by definition.
Since we only consider spectra of $F$-finite rings, $F$-regularity and $F$-purity are preserved under completion.
We omit ``normal'' for $F$-regular singularities because $F$-regularity implies normality \cite{BH}.

The proof of the main theorem heavily depends on Hara's classification of $F$-singularities.
In order to quote results of Hara, we define the ``type'' of a star-shaped weighted dual graph.

\begin{defn}
A \textit{center} of a graph $\Gamma$ is a vertex $v$ having three or more edges directly connected.
$\Gamma$ is a \textit{chain} if it is connected and has neither a center nor a loop.
$\Gamma$ is \textit{star-shaped} if it is connected, has just one center and contains no loop.

If a weighted dual graph $\Gamma$ is star-shaped, we can define the ``\textit{type}'' of $\Gamma$ as follows.
	For each branch $\{v_i\}_{i \in I} (I \subset \{ 1, 2, \cdots n\})$, which is a connected component of $\Gamma$ the unique center removed,
		the \textit{type} of this branch is defined as $d := \mathrm{det} (- E_i \cdot E_j)_{i, j \in I}$.
	If $\Gamma$ has branches of type $d_1, d_2, \cdots d_l \,( d_1 \le d_2 \le \cdots \le d_l)$, $\Gamma$ has \textit{type} $(d_1, d_2 , \cdots d_l)$.
\end{defn}

Some information is omitted compared to Hara's version \cite{Har1}, but it does not matter in almost all cases.

Note that if $(X, x)$ has only a rational singularity, every irreducible component of $E$
	is isomorphic to the projective line $\mathbb{P}_k^1$.
So arbitrary three points on $E_i$ can be taken as $0, - 1, \infty \in \mathbb{P}_k^1$ by an appropriate coordinate change.
On the other hand, four distinct points on $E_i$ can be written as $0, - 1, \lambda, \infty \in \mathbb{P}_k^1\quad(\lambda \neq 0, - 1, \infty)$.

Now we can describe Hara's theorems on $F$-singularities and their dual graphs.

\begin{prop}[{\cite[Theorem (1.1)]{Har1}}]
$(X, x)$ is $F$-regular if and only if it has only a rational singularity and one of the following holds:
\begin{enumerate}
	\item $\Gamma$ is a chain.
	\item $\Gamma$ is star-shaped and either of type
	\begin{enumerate}
		\item $(2, 2, d) (d \ge 2)$, $p \neq 2$,
		\item $(2, 3, 3)$ or $(2, 3, 4)$, $p \neq 2, 3$ or
		\item $(2, 3, 5)$, $p \neq 2, 3, 5$.
	\end{enumerate}
\end{enumerate}
\end{prop}

\begin{figure}[h]
	\centering
	\includegraphics[height=3cm]{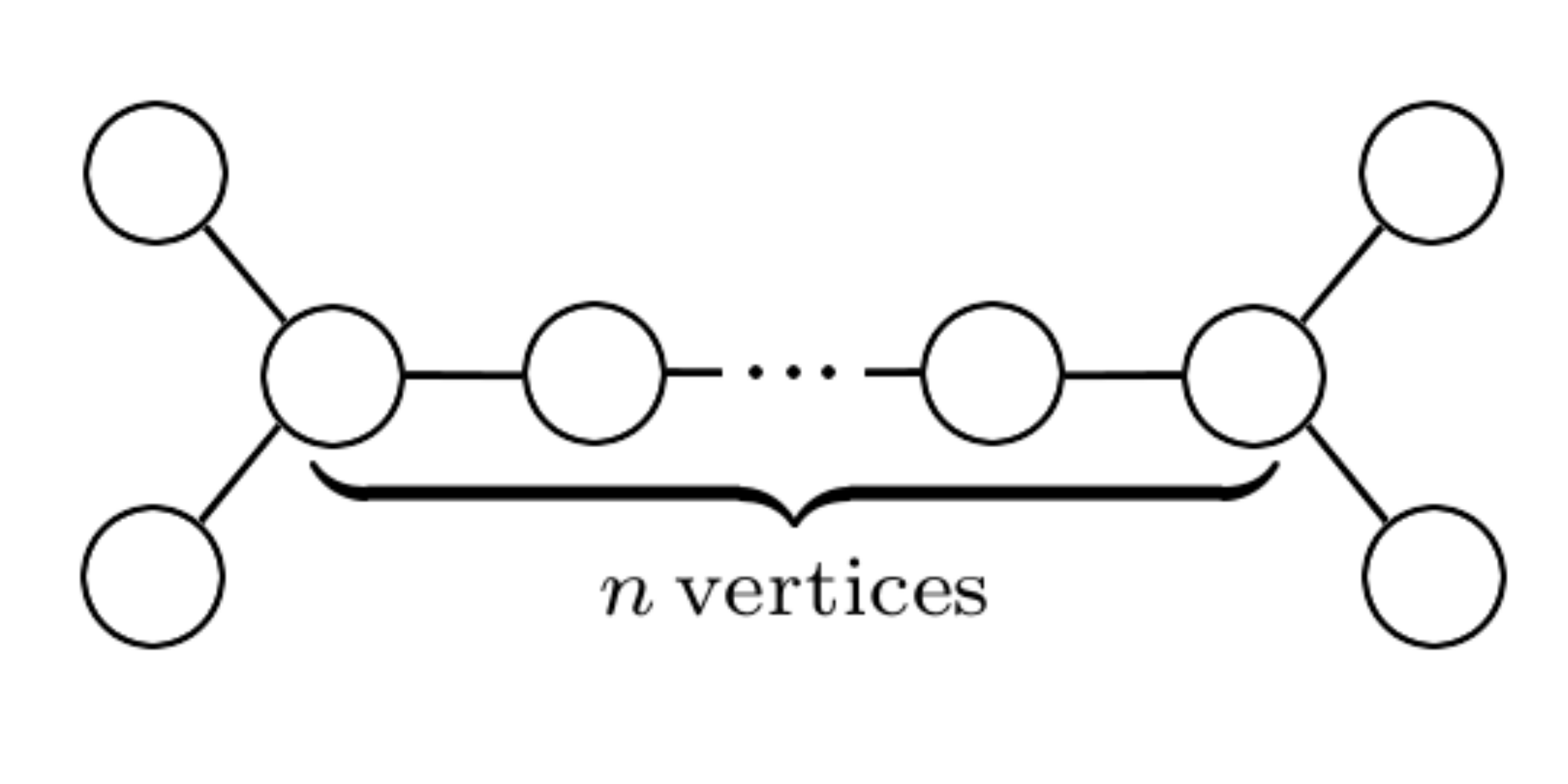}
	\caption{$_*\tilde{D}_{n + 3}$ graph}
	\label{fig:_Dn+3}
\end{figure}

\begin{prop}[{\cite[Theorem (1.2)]{Har1}}]
Assume that $(X, x)$ has only a rational singularity.
If $(X, x)$ is $F$-pure, then one of the following holds:
\begin{enumerate}
	\item $(X, x)$ is $F$-regular
	\item $(X, x)$ is a rational double point, and the graph is either
		\begin{enumerate}
			\item $D_{n+2} (n \ge 2)$, $p = 2$,
			\item $E_6$ or $E_7$, $p = 2, 3$ or
			\item $E_8$, $p = 2, 3, 5$.
		\end{enumerate}
	\item The graph is star-shaped of type either
		\begin{enumerate}
			\item $(3, 3, 3)$ or $(2, 3, 6)$, $p \equiv 1 (\mathrm{mod}\, 3)$,
			\item $(2, 4, 4)$, $p \equiv 1 (\mathrm{mod}\, 4)$ or
			\item $(2, 2, 2, 2)$, $p \neq 2$ and satisfies the condition $(\ast)$. (explained later)
		\end{enumerate}
	\item The graph is $_*\tilde{D}_{n + 3} (n \ge 2)$, $p \neq 2$.(Figure \ref{fig:_Dn+3})
\end{enumerate}
Conversely, if (1), (3) or (4) holds, then $(X, x)$ is $F$-pure.
\end{prop}

Condition $(\ast)$ in (3)(c) is the following:
if we write the intersection points at the central curve as $0, - 1, \lambda, \infty \in \mathbb{P}_k^1$ and $p = 2 m + 1$,
	then the coefficient of $x^m$ in the expansion of $(x + 1)^m (x - \lambda)^m$ is not zero.
Equivalently, $\sum_{k=0}^m \binom{m}{k}^2 (- \lambda)^k \neq 0$ in $k$.
This condition is an open condition for $\lambda \in k$.
In particular, this holds for infinitely many $\lambda$ since $k$ is algebraically closed and therefore infinite field.

If all $E_i$s have the self-intersection $-2$ or less, the type of a branch in a star-shaped graph is strictly larger than the length of the branch.
In other words, the length of a branch is bounded above by its type minus one.
This can be shown by the induction on the length of the branch.
This will help you illustrate the graphs appearing in the above theorems.
For example, (2)(a) case in Proposition 2.3 corresponds to $D_n\, (n \ge 4)$ graphs with the self-intersection $-2$ for length $1$ branches
	and $-2$ or less for the other components.

\section{Tautness criterion}

For two-dimensional normal singularities over $\mathbb{C}$, Laufer gave an equivalent condition to its tautness \cite{Lau1}.
In positive characteristic case, this was partly extended by Sch\"{u}ller \cite{Sch}.
We describe this criterion in this section.

Let $(X, x)$ be a given two-dimensional normal singularity over an algebraically closed field $k$ of positive characteristic $p$
	and $\Gamma_{X, x}$ be the associated weighted dual graph.
Tautness of a nonsingular point is obvious
	and we may assume $\Gamma_{X, x}$ is not a empty graph.

There is a necessary condition for tautness.
\begin{defn}[\cite{Sch}]
	A weighted dual graph $\Gamma$ is \textit{potentially taut} if
		(i) every vertex is associated with the arithmetic genus 0 and
		(ii) every vertex has $3$ or less edges connected directly.
\end{defn}

\begin{prop}[{\cite[Theorem 3.9., Theorem 3.10.]{Lau1}, \cite[Lemma 1.8.]{Sch}}]
	$\Gamma_{X, x}$ is potentially taut if it is taut.
\end{prop}

By this, we may assume that $\Gamma_{X, x}$ is potentially taut and combining this to \cite{LN} not a chain.
Furthermore, we may assume the original singularity is $F$-pure rational in our argument and thus all self-intersection number is at most $- 2$.
We describe properties of $\Gamma_{X, x}$ using not the language of the graphs but of divisors to help you imagine the resulting scheme $P$.
To apply the tautness criterion, we have to construct a \textit{``plumbing scheme''}.

\subsection{Constructing $\Gamma$ from $\Gamma_{X, x}$}
Giving appropriate multiplicities for components of $\Gamma_{X, x}$, we construct a weighted dual graph $\Gamma$.
Since the intersection matrix $\{E_i \cdot E_j\}_{i, j}$ is negative definite \cite{Bad} and in particular invertible matrix,
	there exists an anti-ample cycle $\tilde{Z} = \sum_{i = 1}^n \tilde{\nu_i} E_i \in \mathrm{Div}\, Y$, that is, a cycle satisfying $\tilde{Z} \cdot E_i < 0$ for $i = 1, 2, \cdots n$.
Following the argument on fundamental cycles in \cite{Bad}, $\tilde{\nu_i} > 0$ for all $i$.
In particular, changing $\tilde{Z}$ to its multiple and adding small effective divisor, we may assume $\mathrm{gcd}(\tilde{\nu_i}, p) = 1$.
We fix a sequence of effective divisors
\begin{equation*}
	0 = \tilde{Z}_0 < \tilde{Z}_1 < \cdots < \tilde{Z}_m= \tilde{Z} \quad(\tilde{Z}_{k + 1} = \tilde{Z}_k + E_{i_k}\, (k = 0, 1, \cdots m - 1)).
\end{equation*}

We need some values to construct $\Gamma$.
\begin{defn}
	For an anti-ample divisor $\tilde{Z}$ and a sequence $\{\tilde{Z}_k\}$ as above,
	\begin{equation*}
		\tau := \mathrm{max}_{0 \le k < m}(\tilde{Z}_k \cdot E_{i_k}), \quad\lambda := \mathrm{max}_{1 \le i \le n}(0, 2(2 p_a(E_i) - 2), 2 p_a(E_i) - 2 - E_i^2).
	\end{equation*}
	Since $\tau \ge 1$ and $\lambda \ge 0$ hold in our situations, we assume them.
	Then define the \textit{significant multiplicity} $\nu$ for $\tilde{Z}$ as
	\begin{equation*}
		\nu := \mathrm{min}\{\nu' \in \mathbb{Z} | \nu' \ge \tau + \lambda + 1,\, \mathrm{gcd}(\nu', p) = 1\}
	\end{equation*}
\end{defn}

In \cite{Sch}, more complicated conditions are required for $\nu$.
But these are simplified in our situations.
Note that all coefficients of $\nu \tilde{Z}$ are not divisible by $p$.
Let $Z := \nu \tilde{Z} = \sum_{i = 1}^n \nu_i E_i$ and $\Gamma := \Gamma_Z$.

\subsection{Constructing a plumbing scheme $P$ from $\Gamma$}
Let $\Gamma$ be the weighted dual graph constructed above.
$\Gamma$ has no loop by Proposition 2.4.
For each $E_i$, we construct an open neighborhood $W_i$ of $E_i \subset P$ and glue them into a plumbing scheme $P$.
Then $P$ is a projective scheme embedded in a regular two-dimensional scheme.

First we construct $W_i$.
Let $c_i$ be the number of irreducible components meeting $E_i$.
Since $\Gamma$ is potentially taut and not a chain, $1 \le c_i \le 3$ for all $i$.
Assume $i = 0$ for simplicity.
$W_0$ is defined as a union of two affine schemes
\begin{align*}
	V_0 &= \mathrm{Spec}\, ( k [x_0, y_0, (y_0 - 1)^{- 1}] / (f))\\
	V_0' &= \mathrm{Spec}\, (k [x_0', y_0', (x_0' - 1)^{- 1}, (y_0' - 1)^{- 1}] / (f'))
\end{align*}
where $f \in k [x_0, y_0, (y_0 - 1)^{- 1}]$ and $f' \in k [x_0', y_0', (x_0' - 1)^{- 1}, (y_0' - 1)^{- 1}]$ are polynomials defined below.
Then $V_0$ and $V_0'$ are one-dimensional schemes embedded in $\{y_0 \neq 1\} \subset \mathbb{A}_k^2$
	and $\{x_0' \neq 1, y_0' \neq 1\} \subset \mathbb{A}_k^2$ respectively.
$f$ and $f'$ differ depending on the value of $c_0$.

\begin{enumerate}
	\item If $c_0 =  1$ and $E_1 \cap E_0 \neq \phi$,
		\begin{equation*}
			f = x_0^{\nu_1} y_0^{\nu_0}, \quad f' = {y_0'}^{\nu_0}.
		\end{equation*}
	\item If $c_0 = 2$ and $E_j \cap E_0 \neq \phi \,(j = 1, 2)$,
		\begin{equation*}
			f = x_0^{\nu_1} y_0^{\nu_0}, \quad f' = {x_0'}^{\nu_2} {y_0'}^{\nu_0}.
		\end{equation*}
	\item If $c_0 = 3$ and $E_j \cap E_0 \neq \phi \,(j = 1, 2, 3)$,
		\begin{equation*}
			f = (x_0 - 1)^{\nu_3} x_0^{\nu_1} y_0^{\nu_0}, \quad f' = {x_0'}^{\nu_2} {y_0'}^{\nu_0}.
		\end{equation*}
\end{enumerate}
Then glue them on $\{x_0 \neq 0, 1\} \subset V_0$ and $\{x_0' \neq 0\} \subset V_0'$
	via the coordinate change given by
\begin{equation*}
	\begin{array}{l} 
		x_0 = {x_0'}^{- 1} \\
		y_0 = {x_0'}^{b_0} y_0'
	\end{array}
	\quad \text{and} \quad
	\begin{array}{l} 
		x_0' = x_0^{- 1} \\
		y_0' = x_0^{b_0} y_0
	\end{array}
\end{equation*}
	where $- b_0 = E_0$ is the self-intersection number.
In (3) and (4), there is ambiguity in choice of the order of $E_1, E_2, E_3$.
Although this choice may result a different affine charts, we may choose one arbitrary order.

At the same time, we can glue the neighborhoods of $V_0$ and $V_0'$ into a nonsingular rational surface by the same coordinate change.
So we obtain $W_i$ as a one-dimensional scheme embedded in a nonsingular surface.
$W_i$ has one irreducible component whose reduced structure is isomorphic to $\mathbb{P}_k^1$ and its self-intersection number is $- b_i$.

Now we glue $\{W_i\}_i$ and their neighborhoods into one respectively to obtain $P$.
Assume $\{q\} := E_i \cap E_j \neq \phi (i \neq j)$ and consider the glueing of $W_i$ and $W_j$.
First take a new coordinate system $(x_{i j}, y_{i j})$ on $W_i$ near $q$.
\begin{enumerate}
	\item If $q = \{ x_i = 0 \} \in E_i$, $x_{ij} = x_i, y_{ij} = y_i$ on $\{x_i \neq 1, \infty\} \subset W_i$
	\item If $q = \{ x_i = 1 \} \in E_i$, $x_{ij} = x_i - 1, y_{ij} = y_i$ on $\{x_i \neq 0, \infty\} \subset W_i$
	\item If $q = \{ x_i = \infty \} \in E_i$, $x_{ij} = x_i', y_{ij} = y_i'$ on $\{x_i \neq 0, 1\} \subset W_i$
\end{enumerate}
Take a coordinate system $(x_{j i}, y_{j i})$ on $W_j$ near $q$ in the same way.

Then we can glue appropriate open subsets of them via
\begin{equation*}
	x_{ij} = y_{ji},\quad y_{ij} = x_{ji}.
\end{equation*}
As the construction of $W_i$, we can glue neighborhoods of them at the same time.
Glueing all $W_i$, we obtain $P$ embedded in a regular two-dimensional scheme.
This neighborhood is not necessarily separated.
Now $P$ is a one-dimensional projective scheme over $k$ associated with the weighted dual graph $\Gamma$ as a divisor.

\subsection{Sch\"{u}ller's criterion}

\begin{prop}[{\cite[Proposition 3.16.]{Sch}}]
Let $(X, x)$ be a two-dimensional $F$-pure rational singularity over an algebraically closed field $k$ of positive characteristic $p$.
Assume that the weighted dual graph $\Gamma_{X, x}$ has at least two vertices.
Let $\Gamma$ be the weighted dual graph constructed in Section 3.1. and $P$ be the plumbing scheme for $\Gamma$.
Let $\Theta_P := \mathcal{H}\hspace{-2pt}\mathit{om}_{\mathcal{O}_P}(\Omega_P, \mathcal{O}_P)$ be the tangent sheaf of $P$ where $\Omega_P$ is the sheaf of differentials.
Then $(X, x)$ is taut if $H^1(P, \Theta_P) = 0$.
\end{prop}

\begin{rem}
We can apply similar arguments in the case $k = \mathbb{C}$.
In this case, $H^1(P, \Theta_P) = 0$ is also the necessary condition for $(X, x)$ to be taut \cite{Lau1}.
Sch\"{u}ller's conjecture says that this also holds in positive characteristic.
\end{rem}

We denote $\Gamma(\Theta_U) = \Gamma(U, \Theta_P|_U)$ and $H^1(\Theta_U) = H^1(U, \Theta_P|_U)$ for an open subset $U \subset P$.

\section{Proof of the main theorem}

We prove the main theorem using Sch\"{u}ller's criterion (Proposition 3.4) and Hara's classification (Proposition 2.3).

\subsection{Chain case}
According to \cite{LN}, a rational singularity associated with a chain graph is taut.
This can also be proved by the computation below.

\subsection{Star-shaped case}
The proof requires long computation.

\subsubsection{Forms of each branch}
The following is the list of possible branches of each type.
The number in a vertex represents its self-intersection number.
Self-intersection number $- 2$ is omitted.
\begin{figure}[!h]
	\centering
	\includegraphics[height=1cm]{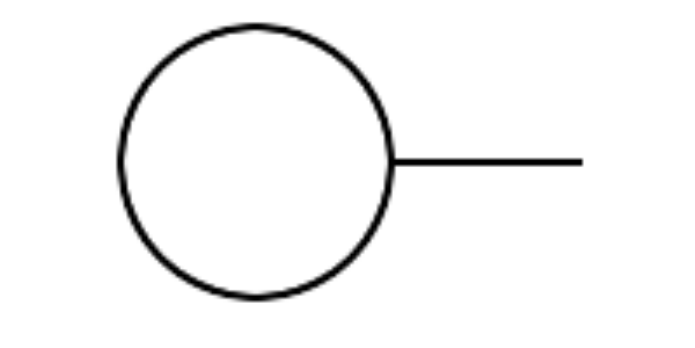}
	\caption{type $2$ branch.}
	\label{type2branch}
\end{figure}%
\begin{figure}[!h]
	\centering
	\includegraphics[height=1cm]{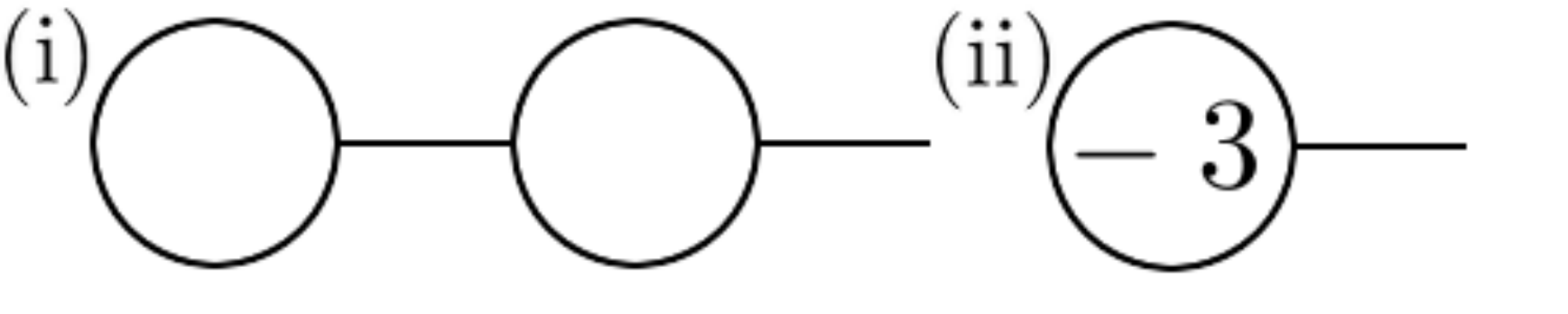}
	\caption{type $3$ branches.}
	\label{type3branch}
	\end{figure}%
\begin{figure}[!h]
	\centering
	\includegraphics[height=1cm]{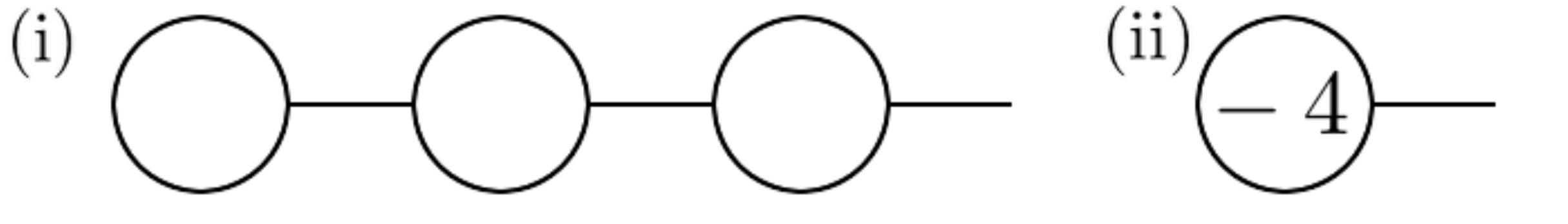}
	\caption{type $4$ branches.}
	\label{type4branch}
\end{figure}%
\begin{figure}[h]
	\centering
	\includegraphics[height=2cm]{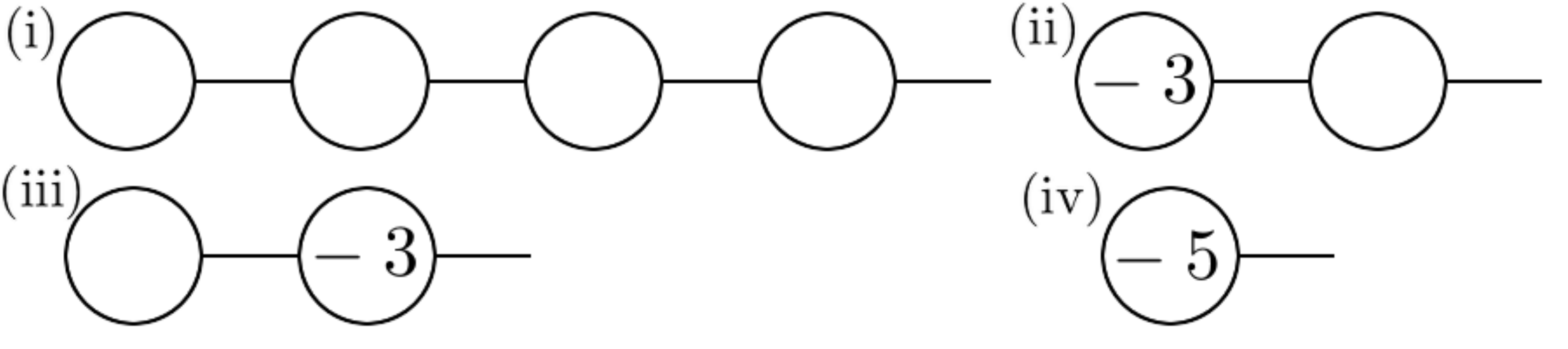}
	\caption{type $5$ branches.}
	\label{type5branch}
\end{figure}%
		
\subsubsection{General settings}

We have to compute $H^1(\Theta_P)$ for all possible cases.
We take an open covering of each plumbing scheme in a common manner.

\begin{figure}[h]
		\centering
		\includegraphics[height=3cm]{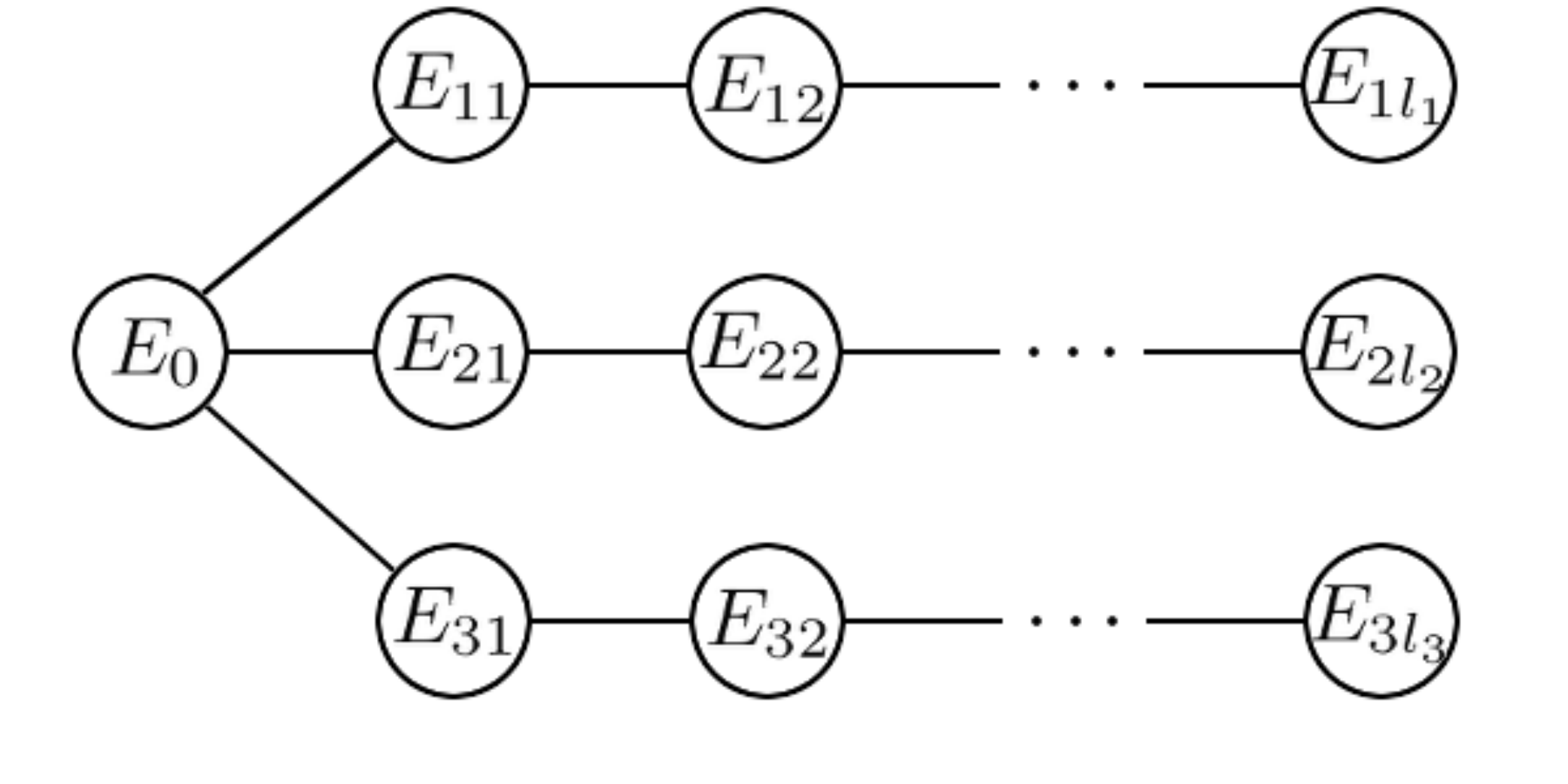}
		\caption{Labeling for a star-shaped graph.}
		\label{fig:starshaped}
	\end{figure}

We fix the notation as follows.
Let the central curve be $E_0$.
Label three branches by $1, 2, 3$ in ascending order of their types.
For branches of the same type, label them in ascending order of the labelings of the branches listed above.
We set the labeling of each irreducible component of P as follows : the component in the $i$-th branch next to $E_0$ is $E_{i 1}$, the next is $E_{i 2}$,
	and the last is $E_{i l_i}$. (Figure \ref{fig:starshaped})
Consequently $l_1 + l_2 + l_3 + 1 = n$.

\begin{figure}[h]
	\centering
	\includegraphics[height=3cm]{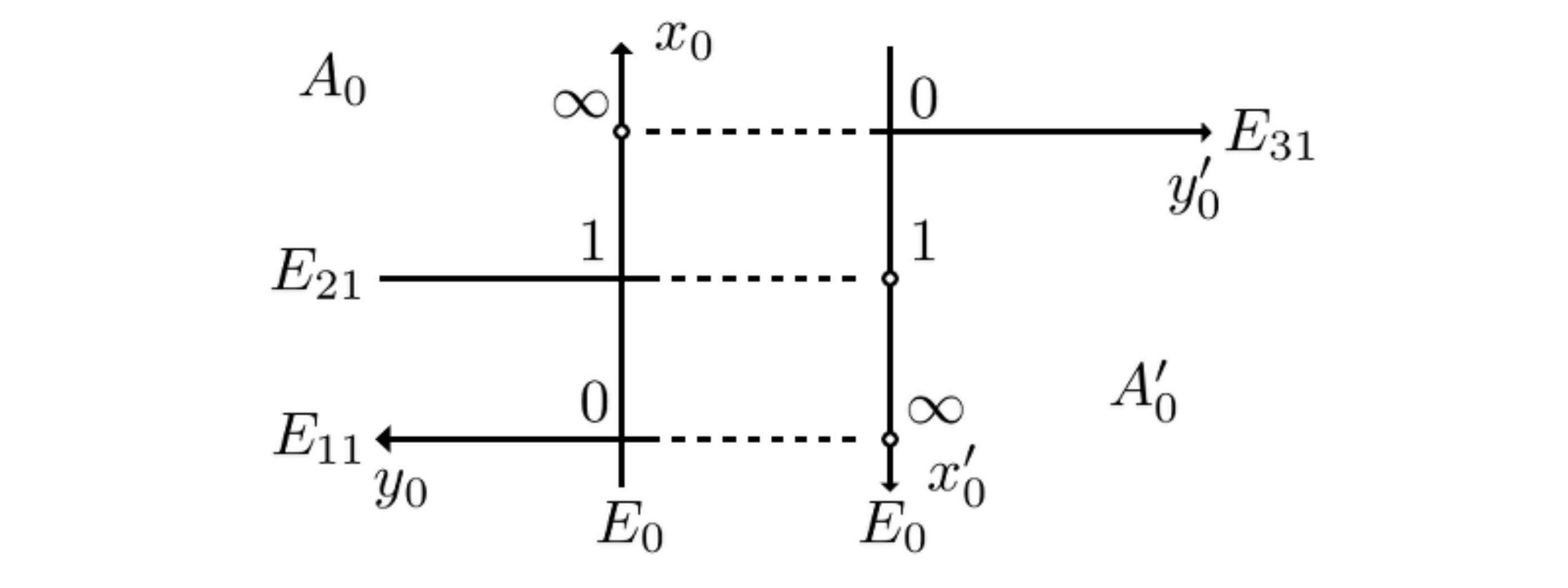}
	\caption{Open affine cover of $E_0 \subset P$.}
	\label{fig:E0}
\end{figure}

Let the intersection of $E_0$ and the first (resp. second, third) branch as $\{x_0 = 0\}$ (resp. $\{x_0 = 1\}, \{x_0 = \infty\}$)
	where $x_0$ is a coordinate of $E_0 \cong \mathbb{P}_k^1$.
Precisely, we cover $E_0 \subset P$ by two open affine subsets $A_0$ and $A_0'$ defined by
\begin{align*}
	A_0 &= \mathrm{Spec}(k[x_0, y_0] / (x_0 - 1)^{\nu_{2 1}} x_0^{\nu_{1 1}} y_0^{\nu_0})\\
	A_0' &= \mathrm{Spec}(k[x_0', y_0', (x_0' - 1)^{- 1}] / {x_0'}^{\nu_{3 1}} {y_0'}^{\nu_0})
\end{align*}
Here $(x_0, y_0)$ and $(x_0', y_0')$ correspond to the coordinates in the construction of $W_0$ (Figure \ref{fig:E0}).

Next we take an open affine covering $\{A_{i j}\}_{j = 1}^{l_i}$ of the $i$-th branch.
For simplicity of notation, we write $E_0$ as $E_{i 0}$.
First let
\begin{equation*}
	A_{i l_i} = E_{i l_i} \setminus (E_{i l_i} \cap E_{i (l_i - 1)}) = \mathrm{Spec}(k[x_{i l_i}, y_{i l_i}] / y_{i l_i}^{\nu_{i l_i}}).
\end{equation*}

We take an open neighborhood of $E_{i (j + 1)} \cap E_{i j} \, (1 \le j \le l_i - 1)$ as
\begin{align*}
	A_{i j} &= (E_{i (j + 1)} \setminus (E_{i (j + 1)} \cap E_{i (j + 2)})) \cup (E_{i j} \setminus (E_{i j} \cap E_{i (j - 1)})) \\
		&= \mathrm{Spec}(k[x_{i j}, y_{i j}] / x_{i j}^{\nu_{i (j + 1)}} y_{i j}^{\nu_{i j}})
\end{align*}
Here $E_{i l_i + 1}$ refers to a point $\{ x_{i l_i} = 0 \} \in E_{i l_i} \setminus E_{i (l_i - 1)}$.
Then $A_{i j}$ is a union of two thickened $\mathbb{A}_k^1$ and $A_{i l_i}$ is a thickened $\mathbb{A}_k^1$.

We have got an open affine covering $\{A_0, A_0'\} \cup \{A_{1 j}\}_j \cup \{A_{2 j}\}_j \cup \{A_{3 j}\}_j$ of $P$.
Let $U_0 = (\bigcup_j A_{1 j}) \cup (\bigcup_j A_{2 j}) \cup A_0$ and $U_1 = \bigcup_j A_{3 j} \cup A_0'$.
Then $U_0 \cap U_1$ is a thickened $\mathbb{P}_k^1$ three points removed.
So $P = U_0 \cup U_1$ is a Leray cover for $\Theta_P$ provided $H^1(\Theta_{U_i}) = 0\, (i = 0, 1)$.
If this is the case,
\begin{equation*}
	H^1(\Theta_P) = \mathrm{Coker}\left(\Gamma(\Theta_{U_0}) \oplus \Gamma(\Theta_{U_1})
		\rightarrow \Gamma(\Theta_{U_0 \cap U_1})\right). \label{eq:H1result}
\end{equation*}

The coordinate changes are given by
\begin{align*}
	x_{i (j + 1)} &= y_{i j}^{- 1}, & y_{i (j + 1)} &= x_{i j} y_{i j}^{b_{i (j + 1)}} \quad(1 \le j \le l_i - 1)\\
	x_{1 1} &= y_0^{- 1}, & y_{1 1} &= x_0 y_0^{b_{1 1}}\\
	x_{2 1} &= y_0^{- 1}, & y_{2 1} &= (x_0 - 1) y_0^{b_{2 1}}\\
	x_{3 1} &= {y_0'}^{- 1}, & y_{3 1} &= {x_0'} {y_0'}^{b_{3 1}}\\
	x_0' &= x_0^{- 1} & y_0' &= x_0^{b_0} y_0.
\end{align*}
This will be used later.

\subsubsection{Local calculation of the tangent sheaf}
Now we start the computation of cohomologies.
First we have to compute the sections of $\Theta_P$ on each affine subsets.
For an affine scheme $A = \mathrm{Spec}(k[x, y] / (f))\, (f \in k[x, y])$, we have an exact sequence
\begin{equation*}
	(f) / (f^2) \stackrel{d}{\rightarrow} (k[x, y] / (f)) dx \oplus (k[x, y] / (f)) dy \rightarrow \Omega_{k[x, y] / (f)} \rightarrow 0 \label{eq:diff}
\end{equation*}
	by the embedding $A \subset \mathbb{A}_k^2$ \cite{Mat}.
Taking $k[x, y]$-dual (or equivalently $k[x, y] / (f)$-dual) of this sequence, we obtain
\begin{equation*}
	0 \rightarrow \Theta_{k[x, y] / (f)} \rightarrow (k[x, y] / (f)) \frac{\partial}{\partial x} \oplus (k[x, y] / (f)) \frac{\partial}{\partial y}. \label{eq:tangent}
\end{equation*}
By this sequence, elements of $\Theta_{k[x, y] / (f)}$ can be represented
	as a $k[x, y] / (f)$-linear sum of $\partial / \partial x$ and $\partial / \partial y$.

\begin{itemize}
	\item If $f = y^{\nu}\, (\mathrm{gcd}(\nu, p) = 1)$, then $df = \nu y^{\nu - 1} dy$ and
		\begin{equation}
			\Theta_{k[x, y] / (f)} = (k[x, y] / (y^{\nu})) \frac{\partial}{\partial x} \oplus (k[x, y] / (y^{\nu - 1})) y \frac{\partial}{\partial y}. \label{theta0}
		\end{equation}
	\item If $f = x^{\nu'} y^{\nu}\, (\mathrm{gcd}(\nu, p) = \mathrm{gcd}(\nu', p) = 1)$, then $df = \nu' x^{\nu' - 1} y^{nu} dx + \nu x^{\nu'} y^{\nu - 1} dy$ and
		\begin{equation}
			\Theta_{k[x, y] / (f)} = (k[x, y] / (x^{\nu' - 1} y^{\nu})) x \frac{\partial}{\partial x} \oplus (k[x, y] / (x^{\nu'} y^{\nu - 1})) y \frac{\partial}{\partial y}. \label{theta1}
		\end{equation}
	\item If $f = (x - 1)^{\nu''} x^{\nu'} y^{\nu}\, (\mathrm{gcd}(\nu, p) = \mathrm{gcd}(\nu', p) = \mathrm{gcd}(\nu'', p) = 1)$,
			then $df = (\nu'' x + \nu' (x - 1)) (x - 1)^{\nu'' - 1} x^{\nu' - 1} y^{nu} dx + \nu (x - 1)^{\nu''} x^{\nu'} y^{\nu - 1} dy$.
		\begin{equation*}
			\left((f) : (\nu'' x + \nu' (x - 1)) (x - 1)^{\nu'' - 1} x^{\nu' - 1} y^{\nu}\right) = ((x - 1) x)
		\end{equation*}
		and we obtain
		\begin{equation}
			\begin{aligned}
			\Theta_{k[x, y] / (f)} = &(k[x, y] / ((x - 1)^{\nu'' - 1} x^{\nu' - 1} y^{\nu})) (x - 1) x \frac{\partial}{\partial x}\\
				&\oplus (k[x, y] / ((x - 1)^{\nu''} x^{\nu'} y^{\nu - 1})) y \frac{\partial}{\partial y}.
			\end{aligned} \label{theta2}
		\end{equation}
\end{itemize}

The coordinate change of differential operators is given as follows:
if $x = {y'}^{- 1}$ and $y = x' {y'}^b$, then $dx = - {y'}^{- 2} dy',\, dy = {y'}^b dx' + b x' {y'}^{b - 1} dy'$ and
\[\frac{\partial}{\partial x} = b x' y' \frac{\partial}{\partial x'} - {y'}^2 \frac{\partial}{\partial y'}, \quad
	\frac{\partial}{\partial y} = {y'}^{- b} \frac{\partial}{\partial x'}.\]
To simplify the notations, linear terms
\[x \frac{\partial}{\partial x} = b x' \frac{\partial}{\partial x'} - y' \frac{\partial}{\partial y'}, \quad
	y \frac{\partial}{\partial y} = x' \frac{\partial}{\partial x'}\]
are convenient.
By $- x \frac{\partial}{\partial x} + b y \frac{\partial}{\partial y} = y' \frac{\partial}{\partial y'}$, we obtain the following lemma.

\begin{lem}
As a $k$-vector space,
\[k \cdot x \frac{\partial}{\partial x} \oplus k \cdot y \frac{\partial}{\partial y}
	= k \cdot x' \frac{\partial}{\partial x'} \oplus k \cdot y' \frac{\partial}{\partial y'}\]
\end{lem}

\subsubsection{Differential forms on three branches}

Now we can calculate the cohomology of the tangent sheaf on each branch.
For simplicity, denote $A_{1 0} = A_{2 0} = A_0$ and $A_{3 0} = A_0'$.

\begin{lem}
$H^1(\Theta_{\bigcup_{j = 1}^{l_i} A_{i j}}) = 0$ and
$\mathrm{Im}(\Gamma(\Theta_{\bigcup_{j = 1}^{l_i} A_{i j}}) \rightarrow \Gamma(\Theta_{A_{i 0} \cap A_{i 1}}))$ has a basis as follows:
\begin{itemize}
	\item If $i = 1$,
		\begin{equation}
			\begin{aligned}
				x_0^s y_0^t x_0 \frac{\partial}{\partial x_0}\quad &(0 \le s \le \nu_{1 1} - 2,\, t \in \mathbb{Z},\, \beta_1 t \le \alpha_1 s)\\
				x_0^s y_0^t y_0 \frac{\partial}{\partial y_0}\quad &(0 \le s \le \nu_{1 1} - 1,\, t \in \mathbb{Z},\, \beta_1 t \le \alpha_1 s)\\
				x_0^s y_0^t (\alpha_1 x_0 \frac{\partial}{\partial x_0} - \beta_1 y_0 \frac{\partial}{\partial y_0})\quad
					&(0 \le s \le \nu_{1 1} - 1,\, \beta_1 t = \alpha_1 s + 1)
			\end{aligned} \label{branch1}
		\end{equation}
	\item If $i = 2$,
		\begin{equation}
			\begin{aligned}
			(x_0 - 1)^r y_0^t (x_0 - 1) \frac{\partial}{\partial x_0}\quad &(0 \le r \le \nu_{2 1} - 2,\, t \in \mathbb{Z},\, \beta_2 t \le \alpha_2 r)\\
			(x_0 - 1)^r y_0^t y_0 \frac{\partial}{\partial y_0}\quad &(0 \le r \le \nu_{2 1} - 1,\, t \in \mathbb{Z},\, \beta_2 t \le \alpha_2 r)\\
			(x_0 - 1)^r y_0^t (\alpha_2 (x_0 - 1) \frac{\partial}{\partial x_0} - \beta_2 y_0 \frac{\partial}{\partial y_0})\quad
					&(0 \le r \le \nu_{2 1} - 1,\, \beta_2 t = \alpha_2 r + 1)
			\end{aligned} \label{branch2}
		\end{equation}
	\item If $i = 3$,
		\begin{equation}
			\begin{aligned}
			{x_0'}^s {y_0'}^t x_0' \frac{\partial}{\partial x_0'}\quad &(0 \le s \le \nu_{3 1} - 2,\, t \in \mathbb{Z},\, \beta_3 t \le \alpha_3 s)\\
			{x_0'}^s {y_0'}^t y_0' \frac{\partial}{\partial y_0'}\quad &(0 \le s \le \nu_{3 1} - 1,\, t \in \mathbb{Z},\, \beta_3 t \le \alpha_3 s)\\
			{x_0'}^s {y_0'}^t (\alpha_3 x_0' \frac{\partial}{\partial x_0'} - \beta_3 y_0' \frac{\partial}{\partial y_0'})\quad
					&(0 \le s \le \nu_{3 1} - 1,\, \beta_3 t = \alpha_3 s + 1)
			\end{aligned} \label{branch3}
		\end{equation}
\end{itemize}
Here $\alpha_i = \mathrm{det}(- E_{i j} \cdot E_{i j'})_{1 \le j, j' \le l_i}$ and $\beta_i = \mathrm{det}(- E_{i j} \cdot E_{i j'})_{2 \le j, j' \le l_i}$.
If $l_i = 1$, set $\beta_i = 1$.
\end{lem}

\textbf{Proof.}
For simplicity, we prove the lemma in the case $i = 1$ and omit the subscript $i$ so that the branch is covered by $\{A_j\}_{j = 1}^l$.

First note $\alpha$ and $\beta$ appeared above are easily calculated as
\begin{equation*}
	\frac{\alpha}{\beta} = b_1 - \frac{1}{b_2 - \frac{1}{b_3 - \cdots \frac{1}{b_{l_i - 1} - \frac{1}{b_{l_i}}}}}\quad
		(\alpha, \beta > 0,\, \mathrm{gcd}(\alpha, \beta) = 1).
\end{equation*}
It is very easy to show by the induction on $l$ using the expansion of the determinant of the intersection matrix.

We prove the lemma by the induction on $l = l_i \ge 1$.
In the case $l = 1$, the branch is $A_1$ and $H^1(\Theta_{A_1}) = 0$ because $A_1$ is affine and $\Theta_P$ is coherent.
By the equation (\ref{theta0}), $\Gamma(\Theta_{A_1})$ has a set of generators consisting of the following elements:
\begin{equation}
	\begin{aligned}
		x_1^s y_1^t x_1 \frac{\partial}{\partial x_1}
			&= x_0^t y_0^{b_1 t - s} \left(b_1 x_0 \frac{\partial}{\partial x_0} - y_0 \frac{\partial}{\partial y_0}\right) &(s, t\ge 0)\\
		x_1^s y_1^t y_1 \frac{\partial}{\partial y_1} &= x_0^t y_0^{b_1 t - s} x_0 \frac{\partial}{\partial x_0} &(s, t \ge 0) \\
		x_1^{- 1} y_1^t x_1 \frac{\partial}{\partial x_1} &= x_0^t y_0^{b_1 t + 1} \left(b_1 x_0 \frac{\partial}{\partial x_0} - y_0 \frac{\partial}{\partial y_0}\right).
	\end{aligned} \label{coordchange}
\end{equation}
On the other hand, $\Gamma(\Theta_{A_0 \cap A_1})$ has a basis as follows:
\begin{align*}
	x_0^s y_0^t x_0 \frac{\partial}{\partial x_0} &(0 \le s \le \nu_1 - 2, t \in \mathbb{Z})\\
	x_0^s y_0^t y_0 \frac{\partial}{\partial y_0} &(0 \le s \le \nu_1 - 1, t \in \mathbb{Z}).
\end{align*}
Using lemma 4.1., $\mathrm{Im}(\Gamma(\Theta_{A_1}) \rightarrow \Gamma(\Theta_{A_0 \cap A_1}))$ has a basis
\begin{equation}
	\begin{aligned}
	x_0^s y_0^t x_0 \frac{\partial}{\partial x_0}\quad &(0 \le s \le \nu_1 - 2, t \le b_1 s)\\
	x_0^s y_0^t y_0 \frac{\partial}{\partial y_0}\quad &(0 \le s \le \nu_1 - 1, t \le b_1 s)\\
	x_0^s y_0^t \left(b_1 x_0 \frac{\partial}{\partial x_0} - y_0 \frac{\partial}{\partial y_0}\right)\quad &(0 \le s \le \nu_1 - 1, t = b_1 s + 1).
	\end{aligned} \label{length1branch}
\end{equation}
In this case, $\alpha = b_1$ and $\beta = 1$.
So (\ref{length1branch}) coincides the set of generators (\ref{branch1}).

Assume $l \ge 2$ and the lemma holds for any smaller $l$.
Let $\gamma_i = \mathrm{det}(- E_j \cdot E_{j'})_{3 \le j, j' \le l_i}$.
$\Gamma(\Theta_{A_1 \cap A_2})$ has a basis as follows:
\begin{align*}
	x_1^s y_1^t x_1 \frac{\partial}{\partial x_1} &(0 \le s \le \nu_2 - 2, t \in \mathbb{Z})\\
	x_1^s y_1^t y_1 \frac{\partial}{\partial y_1} &(0 \le s \le \nu_2 - 1, t \in \mathbb{Z}).
\end{align*}
Applying the induction hypothesis for $\bigcup_{j = 2}^l A_j$, $H^1(\Theta_{\cup_{j = 2}^l A_j}) = 0$ and
	$\mathrm{Im}(\Gamma(\Theta_{\bigcup_{j = 2}^l A_j}) \rightarrow \Gamma(\Theta_{A_1 \cap A_2}))$ has a set of generators as follows:
\begin{align*}
	x_1^s y_1^t x_1 \frac{\partial}{\partial x_1}\quad &(0 \le s \le \nu_2 - 2,\, t \in \mathbb{Z},\, \gamma t \le \beta s)\\
	x_1^s y_1^t y_1 \frac{\partial}{\partial y_1}\quad &(0 \le s \le \nu_2 - 1,\, t \in \mathbb{Z},\, \gamma t \le \beta s)\\
	x_1^s y_1^t (\beta x_1 \frac{\partial}{\partial x_1} - \gamma y_1 \frac{\partial}{\partial y_1})\quad
		&(0 \le s \le \nu_2 - 1,\, t \in \mathbb{Z},\, \gamma t = \beta s + 1).
\end{align*}

If $t \le 0$, $x_1^s y_1^t x_1 \frac{\partial}{\partial x_1}$ and $x_1^s y_1^t y_1 \frac{\partial}{\partial y_1}$ are all contained in
	$\mathrm{Im}(\Gamma(\Theta_{\bigcup_{j = 2}^l A_j}) \rightarrow \Gamma(\Theta_{A_1 \cap A_2}))$ and other cocycles are contained in $\Gamma(\Theta_{A_1})$.
So the restriction map $\Gamma(\Theta_{\cup_{j = 2}^l A_j}) \oplus \Gamma(\Theta_{A_1}) \rightarrow \Gamma(\Theta_{A_1 \cap A_2})$
	is surjective and $H^1(\Theta_{\cup_{j = 1}^l A_j}) = 0$ since $\{\cup_{j = 2}^l A_j, A_1\}$ is a Leray covering of $\cup_{j = 1}^l A_j$ by the induction hypothesis.
$\Gamma(\Theta_{\cup_{j = 1}^l A_j})$ can be computed as
	$\mathrm{Ker}\left(\Gamma(\Theta_{\cup_{j = 2}^l A_j}) \oplus \Gamma(\Theta_{A_1}) \rightarrow \Gamma(\Theta_{A_1 \cap A_2})\right)$.
Since $(\cup_{j = 2}^l A_j) \cap A_0 = \phi$,
\begin{align*}
	&\mathrm{Im}(\Gamma(\Theta_{\bigcup_{j = 1}^l A_j}) \rightarrow \Gamma(\Theta_{A_0 \cap A_1}))\\
	= &\mathrm{Im}\left(\mathrm{Ker}\left(\Gamma(\Theta_{A_1})
		\rightarrow \mathrm{Coker}\left(\Gamma(\Theta_{\bigcup_{j = 2}^l A_j}) \rightarrow \Gamma(\Theta_{A_1 \cap A_2})\right)\right)
		\rightarrow \Gamma(\Theta_{A_0 \cap A_1})\right).
\end{align*}
This looks complicated, but it says that we have to consider coboundaries coming only from $A_1$ to calculate
	$\mathrm{Im}(\Gamma(\Theta_{\bigcup_{j = 1}^l A_j}) \rightarrow \Gamma(\Theta_{A_0 \cap A_1}))$.

$\mathrm{Ker}\left(\Gamma(\Theta_{A_1}) \rightarrow \mathrm{Coker}\left(\Gamma(\Theta_{\bigcup_{j = 2}^l A_j}) \rightarrow \Gamma(\Theta_{A_1 \cap A_2})\right)\right)$
is generated by
\begin{align*}
	x_1^s y_1^t x_1 \frac{\partial}{\partial x_1}\quad &(0 \le s \le \nu_2 - 2,\, t \ge 0,\, \beta s \ge \gamma t)\\
	x_1^s y_1^t x_1 \frac{\partial}{\partial x_1}\quad &(s \ge \nu_2 - 1,\, t \ge 0)\\
	x_1^s y_1^t y_1 \frac{\partial}{\partial y_1}\quad &(0 \le s \le \nu_2 - 1,\, t \ge 0,\, \beta s \ge \gamma t)\\
	x_1^s y_1^t y_1 \frac{\partial}{\partial y_1}\quad &(s \ge \nu_2,\, t \ge 0)\\
	x_1^s y_1^t (\beta x_1 \frac{\partial}{\partial x_1} - \gamma y_1 \frac{\partial}{\partial y_1})\quad
		&(0 \le s \le \nu_{1 2} - 1,\, t \ge 0,\, \beta s + 1 = \gamma t).
\end{align*}
Then we apply the coordinate change (\ref{coordchange}).
Note that $\alpha / \beta = b_1 - \gamma / \beta$.
Then these are
\begin{align*}
	x_0^s y_0^t \left(b_1 x_0 \frac{\partial}{\partial x_0} - y_0 \frac{\partial}{\partial y_0}\right)
		\quad&(s \ge 0, t \le \mathrm{max}\{\frac{\alpha}{\beta} s, b_1 s - \nu_2 + 1\})\\
	x_0^s y_0^t x_0 \frac{\partial}{\partial x_0}
		\quad&(s \ge 0, t \le \mathrm{max}\{\frac{\alpha}{\beta} s, b_1 s - \nu_2\})\\
	x_0^s y_0^t \left(\alpha x_0 \frac{\partial}{\partial x_0} - \beta y_0 \frac{\partial}{\partial y_0}\right)
		\quad&(s \ge 0,\, \alpha s + 1 = \beta t).
\end{align*}
Changing the basis and sending them to $\Gamma(\Theta_{A_0 \cap A_1})$, 
\begin{equation}
	\begin{aligned}
		x_0^s y_0^t x_0 \frac{\partial}{\partial x_0}
			\quad&(0 \le s \le \nu_1 - 2, t \le \mathrm{max}\{\frac{\alpha}{\beta} s, b_1 s - \nu_2\})\\
		x_0^s y_0^t y_0 \frac{\partial}{\partial x_0}
			\quad&(0 \le s \le \nu_1 - 1, t \le \mathrm{max}\{\frac{\alpha}{\beta} s, b_1 s - \nu_2\})\\
		x_0^s y_0^t \left(b_1 x_0 \frac{\partial}{\partial x_0} - y_0 \frac{\partial}{\partial y_0}\right)
			\quad&(0 \le s \le \nu_1 - 1, t = b_1 s - \nu_2 + 1\})\\
		x_0^s y_0^t \left(\alpha x_0 \frac{\partial}{\partial x_0} - \beta y_0 \frac{\partial}{\partial y_0}\right)
			\quad&(0 \le s \le \nu_1 - 1,\, \alpha s + 1 = \beta t).
	\end{aligned} \label{imagewaste}
\end{equation}
Conditions on exponents can be simplified by the following fact.
\begin{equation*}
	s \le \nu_1 \Rightarrow \frac{\alpha}{\beta} s > b_1 s - \nu_2.
\end{equation*}
The resulting inequality is equivalent to $b_1 s - \nu_2 - \frac{\alpha}{\beta} s  = (b_1 - \frac{\alpha}{\beta}) s - \nu_2  = \frac{\gamma}{\beta} s - \nu_2 < 0$.
We show this by the induction on $l \ge 2$.

If $l = 2$, then $\alpha = b_1 - 1 / b_2,\, \beta = b_2$ and $\gamma = 1$.
\begin{align*}
	s - b_2 \nu_2 &\le \nu_1 - b_2 \nu_2\\
		&= E_2 \cdot \nu_1 E_1 + E_2 \cdot \nu_2 E_2\\
		&= E_2 \cdot \nu \tilde{Z} - E_2 \cdot \nu_0 E_0\\
		&< E_2 \cdot \nu \tilde{Z}\\
		&< 0
\end{align*}
because $\nu \tilde{Z}$ is anti-ample.
So $\frac{\gamma}{\beta} s - \nu_2 = \frac{1}{b_2} (s - b_2 \nu_2) < 0$.

Next assume $l \ge 3$ and this holds for $l - 1$.
Let $\delta$ be $\mathrm{det}(- E_{i j} \cdot E_{i j'})_{4 \le j, j' \le l_i}$ for $l \ge 4$ and $1$ for $l = 3$.
Then $\beta = b_2 \gamma - \delta$ and
\begin{align*}
	\gamma s - \beta \nu_2 &\le \gamma \nu_1 - (b_2 \gamma - \delta) \nu_2\\
		&= \gamma (\nu_1 - b_2 \nu_2) + \delta \nu_2\\
		&= \gamma (E_2 \cdot \nu_1 E_1 + E_2 \cdot \nu_2 E_2) + \delta \nu_2\\
		&= \gamma (E_2 \cdot \nu \tilde{Z} - E_2 \cdot \nu_3 E_3) + \delta \nu_2\\
		&< - \gamma \nu_3 + \delta \nu_2\\
		&< 0
\end{align*}
by the induction hypothesis.

Consequently, coboundaries (\ref{imagewaste}) can be written as
\begin{align*}
	x_0^s y_0^t x_0 \frac{\partial}{\partial x_0}\quad &(0 \le s \le \nu_1 - 2,\, t \in \mathbb{Z},\, \beta t \le \alpha s)\\
	x_0^s y_0^t y_0 \frac{\partial}{\partial y_0}\quad &(0 \le s \le \nu_1 - 1,\, t \in \mathbb{Z},\, \beta t \le \alpha s)\\
	x_0^s y_0^t (\alpha_1 x_0 \frac{\partial}{\partial x_0} - \beta_1 y_0 \frac{\partial}{\partial y_0})\quad
		&(0 \le s \le \nu_1 - 1,\, \beta t = \alpha s + 1),
\end{align*}
which is the desired form. $\square$

\vspace{1\baselineskip}
We need a further coordinate change and computation for coboundaries from $U_1$.
$\Gamma(\Theta_{A_0'})$ has a set of generators consisting of following elements:
\begin{align*}
	{x_0'}^s {y_0'}^t x_0' \frac{\partial}{\partial x_0'} \quad&(s, t \ge 0)\\
	(1 - x_0')^r {y_0'}^t x_0' \frac{\partial}{\partial x_0'} \quad&(r < 0,\, t \ge 0)\\
	{x_0'}^s {y_0'}^t y_0' \frac{\partial}{\partial y_0'} \quad&(s, t \ge 0)\\
	(1 - x_0')^r {y_0'}^t y_0' \frac{\partial}{\partial y_0'} \quad&(r < 0,\, t \ge 0)
\end{align*}
Same computation as the proof of lemma 4.2 shows that the restriction map
\begin{equation}
	\Gamma(\Theta_{A_0'}) \rightarrow \mathrm{Coker}\left(\Gamma(\Theta_{\cup_j A_{3 j}}) \rightarrow \Gamma(\Theta_{A_0' \cap A_{3 1}})\right) \label{restinU1}
\end{equation}
is surjective.
This shows $H^1(\Theta_{U_1}) = 0$ because $U_1 = A_0' \cup (\bigcup_j A_{3 j})$ is a Leray covering for $\Theta_{U_1}$.
The kernel of (\ref{restinU1}) has a basis consisting of
\begin{align*}
	{x_0'}^s {y_0'}^t x_0' \frac{\partial}{\partial x_0'}\quad (t \ge 0,\, s &\ge \mathrm{min}\{\frac{\beta_3}{\alpha_3} t, \nu_{3 1} - 1\})\\
	{x_0'}^s {y_0'}^t y_0' \frac{\partial}{\partial y_0'}\quad (t \ge 0,\, s &\ge \mathrm{min}\{\frac{\beta_3}{\alpha_3} t, \nu_{3 1}\})\\
	{x_0'}^s {y_0'}^t (\alpha_3 x_0' \frac{\partial}{\partial x_0'} - \beta_3 y_0' \frac{\partial}{\partial y_0'})\quad
		(0 \le s \le &\nu_{3 1} - 1,\, \beta_3 t = \alpha_3 s + 1)
\end{align*}
\begin{align*}
	\left((1 - x_0')^r - (1 + x_0' + \cdots + {x_0'}^{\nu_{3 1} - 1})^{- r}\right) {y_0'}^t x_0' \frac{\partial}{\partial x_0'}\quad &(r < 0,\, t \ge 0)\\
	\left((1 - x_0')^r - (1 + x_0' + \cdots + {x_0'}^{\nu_{3 1} - 1})^{- r}\right) {y_0'}^t y_0' \frac{\partial}{\partial y_0'}\quad &(r < 0,\, t \ge 0).
\end{align*}
This is a part of a basis of $\Gamma(\Theta_{U_1})$ containing all elements necessary for computing the coboundaries.
Sending them to $\Gamma(\Theta_{A_0 \cap A_0'})$, these are
\begin{align}
	{x_0'}^s {y_0'}^t x_0' \frac{\partial}{\partial x_0'}\quad &(0 \le t \le \nu_0 - 1,\, \beta_3 t \le \alpha_3 s) \label{term1}\\
	{x_0'}^s {y_0'}^t y_0' \frac{\partial}{\partial y_0'}\quad &(0 \le t \le \nu_0 - 2,\, \beta_3 t \le \alpha_3 s) \label{term2}\\
	{x_0'}^s {y_0'}^t (\alpha_3 x_0' \frac{\partial}{\partial x_0'} - \beta_3 y_0' \frac{\partial}{\partial y_0'})\quad
		&(s \ge 0,\, \beta_3 t = \alpha_3 s + 1) \label{term3}\\
	\left((1 - x_0')^r - (1 + x_0' + \cdots + {x_0'}^{\nu_{3 1} - 1})^{- r}\right) {y_0'}^t x_0' \frac{\partial}{\partial x_0'}\quad &(r < 0,\, 0 \le t \le \nu_0 - 1)\\
	\left((1 - x_0')^r - (1 + x_0' + \cdots + {x_0'}^{\nu_{3 1} - 1})^{- r}\right) {y_0'}^t y_0' \frac{\partial}{\partial y_0'}\quad &(r < 0,\, 0 \le t \le \nu_0 - 2).
\end{align}
This is a basis of $\mathrm{Im}(\Gamma(\Theta_{U_1}) \rightarrow \Gamma(\Theta_{U_0 \cap U_1}))$, coboundaries from $U_1$.

\subsubsection{Differential forms near the central curve.}

In 4.2.4, $H^1(U_1) = 0$ is shown and $\mathrm{Im}(\Gamma(\Theta_{U_1}) \rightarrow \Gamma(\Theta_{U_0 \cap U_1}))$ is computed.
Next we show $H^1(U_0) = 0$ and compute $\mathrm{Im}(\Gamma(\Theta_{U_0}) \rightarrow \Gamma(\Theta_{U_0 \cap U_1}))$.

$H^1(\Theta_{U_0})$ can be calculated as
\[H^1(\Theta_{U_0}) = \mathrm{Coker}\left(\Gamma(\Theta_{A_0})
		\rightarrow \bigoplus_{i = 1, 2} \mathrm{Coker}(\Gamma(\Theta_{\cup_{j = 1}^{l_i} A_{i j}}) \rightarrow \Gamma(\Theta_{A_0 \cap A_{i 1}}))\right)\]
because $U_0 = A_0 \cup (\bigcup_j A_{1 j}) \cup (\bigcup_j A_{2 j})$ is a Leray covering for $\Theta_{U_0}$.
Let
\begin{align*}
	C_1 &:= \mathrm{Coker}(\Gamma(\Theta_{\cup_j A_{1 j}}) \rightarrow \Gamma(\Theta_{A_0 \cap A_{1 1}}))\\
	C_2 &:= \mathrm{Coker}(\Gamma(\Theta_{\cup_j A_{2 j}}) \rightarrow \Gamma(\Theta_{A_0 \cap A_{2 1}}))\\
	K_1 &:= \mathrm{Ker}\left(\Gamma(\Theta_{A_0}) \rightarrow C_1\right) \subset \Gamma(\Theta_{A_0})\\
	K_2 &:= \mathrm{Ker}\left(\Gamma(\Theta_{A_0}) \rightarrow C_2\right) \subset \Gamma(\Theta_{A_0}).
\end{align*}
Then $H^1(\Theta_{U_0})$ vanishes if and only if $K_2 \rightarrow C_1$ and $K_1 \rightarrow C_2$ are both surjective.
Now we prove only surjectivity of the map $K_2 \rightarrow C_1$.
The proof of the other goes symmetrically.
By (\ref{theta2}) and (\ref{branch2}), $K_2$ is generated by the following elements:
\begin{align*}
	(x_0 - 1)^r y_0^t x_0 \frac{\partial}{\partial x_0}\quad &(t \ge 0, r \ge \mathrm{min}(\frac{\beta_2}{\alpha_2} t + 1, \nu_{2 1}))\\
	(x_0 - 1)^r y_0^t y_0 \frac{\partial}{\partial y_0}\quad &(t \ge 0, r \ge \mathrm{min}(\frac{\beta_2}{\alpha_2} t, \nu_{2 1}))\\
	(x_0 - 1)^r y_0^t (\alpha_2 ((x_0 - 1) \frac{\partial}{\partial x_0} - \beta_2 y_0 \frac{\partial}{\partial y_0})\quad
		&(1 \le r \le \nu_{2 1} - 1,\, \beta_2 t = \alpha_2 r + 1)\\
\end{align*}
Changing the basis, this is same as
\begin{equation}
	\begin{aligned}
		(x_0 - 1)^r x_0^s y_0^t x_0 \frac{\partial}{\partial x_0}\quad &(s, t \ge 0, r = \mathrm{min}(\lceil\frac{\beta_2}{\alpha_2} t + 1\rceil, \nu_{2 1}))\\
		(x_0 - 1)^r x_0^s y_0^t y_0 \frac{\partial}{\partial y_0}\quad &(s, t \ge 0, r = \mathrm{min}(\lceil\frac{\beta_2}{\alpha_2} t\rceil, \nu_{2 1}))\\
		\alpha_2 (x_0 - 1)^{r + 1} y_0^t x_0 \frac{\partial}{\partial x_0} - \beta_2 (x_0 - 1)^r y_0^t y_0 \frac{\partial}{\partial y_0}\quad
				&(1 \le r \le \nu_{2 1} - 1,\, \beta_2 t = \alpha_2 r + 1).
	\end{aligned} \label{K2basis}
\end{equation}

We check
\begin{equation}
x_0^s y_0^t x_0 \frac{\partial}{\partial x_0}
	\in \mathrm{Im}\left(\Gamma(\Theta_{\cup_j A_{1 j}}) \oplus K_2 \rightarrow \Gamma(\Theta_{A_0 \cap A_1})\right) \label{inimage}
\end{equation}
for $0 \le s \le \nu_{1 1} - 1$ and $t \in \mathbb{Z}$.
If $s \ge \frac{\beta_1}{\alpha_1} t$, it is already in
	$\mathrm{Im}\left(\Gamma(\Theta_{\cup_j A_{1 j}}) \rightarrow \Gamma(\Theta_{A_0 \cap A_{1 1}})\right)$
by (\ref{branch1}).
We prove (\ref{inimage}) by the descending induction on $s < \frac{\beta_1}{\alpha_1} t$.
Assume (\ref{inimage}) holds for all larger $s$.
If $r = \mathrm{min}(\lceil\frac{\beta_2}{\alpha_2} t + 1\rceil, \nu_{2 1})$,
then $(- 1)^r (x_0 - 1)^r x_0^s y_0^t x_0 \frac{\partial}{\partial x_0} \in \mathrm{Im}(K_2 \rightarrow \Gamma(\Theta_{A_0} \cap A_{1 1}))$
is a sum of $x_0^s y_0^t x_0 \frac{\partial}{\partial x_0}$ and terms with higher degree in $x_0$.
Since these higher degree terms are in the image by the induction hypothesis, $x_0^s y_0^t x_0 \frac{\partial}{\partial x_0}$ is also in the image.
Then (\ref{inimage}) is proved.
Similar argument can be applied to $y_0 \frac{\partial}{\partial y_0}$ terms and surjectivity of $K_2 \rightarrow C_1$ is proved.
Then we obtain $H^1(\Theta_{U_0}) = 0$.

Next we need a basis of
\[\Gamma(\Theta_{U_0}) = \mathrm{Ker}(\Gamma(\Theta_{A_0}) \rightarrow C_1 \oplus C_2) = K_1 \cap K_2\]
but this can be easily obtained by (\ref{K2basis}).
So we give the list later.

\subsubsection{Remarks on $\alpha_i$ and $\beta_i$}

The coboundary space from $U_0$ and $U_1$ are determined by $\alpha_i$ and $\beta_i$ for $i = 1, 2, 3$.
If $0 \le \gamma_i < \beta_i$,
\[\alpha_i = b_{i 1} \beta_i - \gamma_i \ge 2 \beta_i - \gamma_i > \beta_i.\]
Applying same argument for subgraphs of each branch, we obtain $\alpha_i > \beta_i > 0\quad (i = 1, 2, 3)$.
In particular $\beta_i / \alpha_i \le 1 - \frac{1}{d_i}$ because $\alpha_i$ coincides with the type of the $i$-th branch $d_i$.

\subsubsection{Computing $H^1(\Theta_P) = 0$}

Now we show
\[H^1(\Theta_P) = \mathrm{Coker}(\Gamma(\Theta_{U_0}) \oplus \Gamma(\Theta_{U_1}) \rightarrow \Gamma(\Theta_{A_0 \cap A_1})) = 0.\]

We take a basis of $\Gamma(\Theta_{A_0 \cap A_1})$ as follows:
\begin{align}
	x_0^s y_0^t x_0 \frac{\partial}{\partial x_0}\quad &(s \in \mathbb{Z},\, 0 \le t \le \nu_0 - 1) \label{xterm}\\
	x_0^s y_0^t y_0 \frac{\partial}{\partial x_0}\quad &(s \in \mathbb{Z},\, 0 \le t \le \nu_0 - 2) \label{yterm}\\
	(x_0' - 1)^r {y_0'}^t x_0' \frac{\partial}{\partial x_0'}\quad &(r < 0,\, 0 \le t \le \nu_0 - 1) \label{poleterm1}\\
	(x_0' - 1)^r {y_0'}^t y_0' \frac{\partial}{\partial y_0'}\quad &(r < 0,\, 0 \le t \le \nu_0 - 2). \label{poleterm2}
\end{align}

Now $\Gamma(\Theta_{U_0}) = Z_1 \cap Z_2$ has a set of generators
\begin{align}
	&(x_0 - 1)^r x_0^s y_0^t x_0 \frac{\partial}{\partial x_0} \label{U0xterm}\\
		&\left(t \ge 0,\, r = \mathrm{min}\left\{\lceil\frac{\beta_2}{\alpha_2} t + 1\rceil, \nu_{2 1}\right\},\,
			s \ge \mathrm{min}\left\{\frac{\beta_1}{\alpha_1} t, \nu_{1 1} - 1\right\}\right)\notag\\
	&(x_0 - 1)^r x_0^s y_0^t y_0 \frac{\partial}{\partial y_0} \label{U0yterm}\\
		&\left(t \ge 0,\, r = \mathrm{min}\left\{\lceil\frac{\beta_2}{\alpha_2} t\rceil, \nu_{2 1}\right\},\,
			s \ge \mathrm{min}\left\{\frac{\beta_1}{\alpha_1} t, \nu_{1 1}\right\}\right)\notag\\
	&(x_0 - 1)^r x_0^s y_0^t (\alpha_1 x_0 \frac{\partial}{\partial x_0} - \beta_1 y_0 \frac{\partial}{\partial y_0}) \label{branch1rel}\\
			&\left(r, t \ge 1,\, r = \mathrm{min}\left\{\lceil\frac{\beta_2}{\alpha_2} t + 1\rceil, \nu_{2 1}\right\},\, s = \frac{\beta_1}{\alpha_1} t - \frac{1}{\alpha_1}\right)\notag\\
	&\alpha_2 (x_0 - 1)^{r + 1} x_0^s y_0^t x_0 \frac{\partial}{\partial x_0} - \beta_2 (x_0 - 1)^r x_0^s y_0 \frac{\partial}{\partial y_0} \label{branch2rel}\\
			&\left(r, t \ge 1,\, r = \frac{\beta_2}{\alpha_2} t - \frac{1}{\alpha_2},\, s = \mathrm{min}\left\{\lceil\frac{\beta_1}{\alpha_1} t\rceil, \nu_{1 1}\right\}\right)\notag
\end{align}
by (\ref{branch1}) and (\ref{K2basis}).

On the other hand, by the coordinate change of (\ref{term1}), (\ref{term2}) and (\ref{term3}), $\Gamma(\Theta_{U_1})$ has a basis
\begin{align}
	x_0^s y_0^t x_0 \frac{\partial}{\partial x_0}\quad
		&\left(0 \le t \le \nu_0 - 1,\, s \le \frac{\alpha'}{\alpha_3} t\right) \label{U1xterm}\\
	x_0^s y_0^t y_0 \frac{\partial}{\partial y_0}\quad
		&\left(0 \le t \le \nu_0 - 2,\, s \le \frac{\alpha'}{\alpha_3} t\right) \label{U1yterm}\\
	x_0^s y_0^t \left(\alpha_3 x_0 \frac{\partial}{\partial x_0} - \alpha' y_0 \frac{\partial}{\partial y_0}\right)\quad
		&\left(0 \le t \le \nu_0 - 1,\, s = \frac{\alpha'}{\alpha_3} t + \frac{1}{\alpha_3}\right) \label{U1rel}\\
	\left((1 - x_0')^r - (1 + x_0' + \cdots + {x_0'}^{\nu_{3 1} - 2})^{- r}\right) &{y_0'}^t x_0' \frac{\partial}{\partial x_0'}\quad (r < 0,\, 0 \le t \le \nu_0 - 1) \label{polecoboundaryx}\\
	\left((1 - x_0')^r - (1 + x_0' + \cdots + {x_0'}^{\nu_{3 1} - 1})^{- r}\right) &{y_0'}^t y_0' \frac{\partial}{\partial y_0'}\quad (r < 0,\, 0 \le t \le \nu_0 - 2). \label{polecoboundaryy}
\end{align}
where $\alpha' = b_0 \alpha_3 - \beta_3$.

\vspace{1\baselineskip}
First we prove that cocycles of the form $x_0^s y_0^t y_0 \frac{\partial}{\partial y_0}$ are all coboundaries.
If $s \le \frac{\alpha_3}{\alpha'} t$, it is a coboundary coming from $U_1$ by (\ref{U1yterm}).
For $s > \frac{\alpha'}{\alpha_3} t$, we show that $x_0^s y_0^t y_0 \frac{\partial}{\partial y_0}$ is a coboundary by the induction on $s$ fixing $t$.
We show that $(x_0 - 1)^r x_0^{s - r} y_0^t y_0 \frac{\partial}{\partial y_0}$ is a coboundary from $U_0$ where $r = \lceil\frac{\beta_2}{\alpha_2} t\rceil$ later.
If this was shown, $(x_0 - 1)^r x_0^{s - r} y_0^t y_0 \frac{\partial}{\partial y_0}$ is a sum of
	$x_0^s y_0^t y_0 \frac{\partial}{\partial y_0}$ and terms with lower degrees in $x_0$.
By the induction hypothesis, these accompanying terms are coboundaries.
Canceling them, we know that $x_0^s y_0^t y_0 \frac{\partial}{\partial y_0}$ is also a coboundary.

What we have to show is $s - r \ge \frac{\beta_1}{\alpha_1} t$.
It is enough to show
\begin{equation}
	\lfloor\frac{\alpha'}{\alpha_3} t\rfloor + 1
		\ge \frac{\beta_1}{\alpha_1} t + \lceil\frac{\beta_2}{\alpha_2} t\rceil. \label{ytermineq}
\end{equation}
We check this case by case.
\begin{enumerate}
	\item Type $(2, 2, d)$ case:
		since type $2$ branch consists of a $(- 2)$-curve,
			$\frac{\beta_1}{\alpha_1} = \frac{\beta_2}{\alpha_2} = \frac{1}{2}$ and $\frac{\alpha'}{\alpha_3} > 1$.
		So
		\begin{equation*}
			\lfloor\frac{\alpha'}{\alpha_3} t\rfloor + 1
				> t + 1
				> (\frac{1}{2} + \frac{1}{2}) t + \frac{1}{2}
				\ge \frac{\beta_1}{\alpha_1} t + \lceil\frac{\beta_2}{\alpha_2} t\rceil.
		\end{equation*}

	\item Type $(2, 3, 3), (2, 3, 4), (2, 3, 5)$ case:
		$\frac{\beta_1}{\alpha_1} t = \frac{1}{2},\, \frac{\beta_2}{\alpha_2} \le \frac{2}{3}$ and $\frac{\alpha'}{\alpha_3} \ge \frac{6}{5}$.
		By
		\begin{equation*}
			\lfloor\frac{\alpha'}{\alpha_3} t\rfloor + 1
				\ge \frac{6}{5} t + \frac{1}{5}
				\ge \frac{1}{2} t + \lceil\frac{2}{3} t\rceil + \frac{1}{30} t - \frac{7}{15},
		\end{equation*}
		(\ref{ytermineq}) holds if $t \ge 14$.
		Direct calculation shows that (\ref{ytermineq}) also holds for $0 \le t \le 13$.
\end{enumerate}

\vspace{1\baselineskip}
Next we show that cocycles of the form $x_0^s y_0^t x_0 \frac{\partial}{\partial x_0}$ are all coboundaries.
This is much harder than $x_0^s y_0^t y_0 \frac{\partial}{\partial y_0}$ terms and characteristic conditions effects critically.
We have coboundaries of the form (\ref{U0xterm}) and call them type A coboundaries.
In any cases of $F$-regular singularities, $\alpha_i = d_i\, (i = 1, 2)$ is smaller than $p$ and nonzero in $k$.
Since monomial terms $x_0^s y_0^t y_0 \frac{\partial}{\partial y_0}$ are all coboundaries,
	we get following coboundaries by subtracting them from (\ref{branch1rel}) and (\ref{branch2rel}):
\begin{align}
	(x_0 - 1)^r x_0^s y_0^t x_0 \frac{\partial}{\partial x_0}\quad
		&\left(t \ge 1,\, r = \mathrm{min}\left\{\lceil\frac{\beta_2}{\alpha_2} t + 1\rceil, \nu_{2 1}\right\},\,
			s = \frac{\beta_1}{\alpha_1} t - \frac{1}{\alpha_1}\right) \label{branch1rel2}\\
	(x_0 - 1)^r x_0^s y_0^t x_0 \frac{\partial}{\partial x_0}\quad
		&\left(t \ge 1,\, r = \frac{\beta_2}{\alpha_2} t + \frac{\alpha_2 - 1}{\alpha_2},\,
			s = \mathrm{min}\left\{\lceil\frac{\beta_1}{\alpha_1} t\rceil, \nu_{1 1}\right\}\right). \label{branch2rel2}
\end{align}
We call (\ref{branch1rel2}) and (\ref{branch2rel2}) coboundaries of type B and C respectively.
Similarly subtracting $x_0^s y_0^t y_0 \frac{\partial}{\partial y_0}$ terms from (\ref{U1rel}), we get the following coboundaries of type D except in the type $(2, 2, d)$ cases:
\begin{equation*}
	x_0^s y_0^t x_0 \frac{\partial}{\partial x_0}\quad
		\left(0 \le t \le \nu_0 - 1,\, s = \frac{\alpha'}{\alpha_3} t + \frac{1}{\alpha_3}\right).
\end{equation*}
The proof uses the basically same method as $x_0^s y_0^t y_0 \frac{\partial}{\partial y_0}$ part.

If $s \le \frac{\alpha'}{\alpha_3} t$, then it is a coboundary from $U_1$ by (\ref{U1xterm}), in other words type A coboundary.
For $s > \frac{\alpha'}{\alpha_3} t$, we use the induction on $s$.
It is enough to show that there exists an integer $r$ such that
\begin{equation}
	(x_0 - 1)^r x_0^{s - r} y_0^t x_0 \frac{\partial}{\partial x_0}\, \text{is a coboundary.} \label{judge}
\end{equation}
Again we calculate this case by case.
We only consider minimum $\frac{\beta_1}{\alpha_1}, \frac{\beta_2}{\alpha_2}$ and $\frac{\alpha'}{\alpha_3}$ because the coboundary space becomes smallest.
\begin{enumerate}
	\item Type $(2, 2, d)\quad (d \ge 2,\, p \neq 2)$ case:
		then $\frac{\beta_1}{\alpha_1} = \frac{\beta_2}{\alpha_2} = \frac{1}{2}$ and $\frac{\alpha'}{\alpha_3} > 1$.
		If $t$ is even, set $r = \frac{1}{2} t + 1$.
		Then $s - r \ge t + 1 - (\frac{1}{2} t + 1) = \frac{1}{2} t$ and this shows (\ref{judge}) holds by a type A coboundary.
		
		If $t$ is odd, set $r = \frac{1}{2} t + \frac{3}{2}$.
		Then $s - r \ge t + 1 - (\frac{1}{2} t + \frac{3}{2}) = \frac{1}{2} t - \frac{1}{2}$.
		If $s - r \ge \frac{1}{2} t + \frac{1}{2}$, it is a type A coboundary.
		Otherwise, $s - r = \frac{1}{2} t - \frac{1}{2}$ and it is a type B coboundary.
		
	\item Type $(2, 3, 3), (2, 3, 4)\quad (p \neq 2, 3)$ case:
		then $\frac{\beta_1}{\alpha_1} = \frac{1}{2}, \frac{\beta_2}{\alpha_2} \le \frac{2}{3}$ and $\frac{\alpha'}{\alpha_3} \ge \frac{5}{4}$.
		If $r = \lceil\frac{\beta_2}{\alpha_2} t + 1\rceil$,
			$s - r  - \frac{1}{2} t \ge \frac{5}{4} t + \frac{1}{4} - (\frac{2}{3} t + \frac{5}{3}) - \frac{1}{2} t = \frac{1}{12} t - \frac{17}{12}$.
		So (\ref{judge}) holds by a type A coboundary if $t \ge 17$.
		(\ref{judge}) also holds for $0 \le t \le 16$ cases by Table \ref{case2table}.
		\begin{table}[h]
		\caption{Types of coboundaries in case (2).}
		\label{case2table}
		\begin{tabular}{@{}|c|c|c|c|c|c|c|c|c|c|c|c|c|c|c|c|c|c|}\hline
			$t$					&0	&1	&2	&3	&4	&5	&6	&7	&8	&9	&10	&11	&12	&13	&14	&15	&16\\ \hline
			$\mathrm{min}\{s\}$	&1	&2	&3	&4	&6	&7	&8	&9	&11	&12	&13	&14	&16	&17	&18	&19	&21\\ \hline	
			$r$					&1	&2	&2	&3	&4	&5	&5	&6	&7	&7	&8	&9	&9	&10	&11	&11	&12\\ \hline
			$\mathrm{min}\{s - r\}$	&0	&0	&0	&1	&2	&2	&3	&3	&4	&5	&5	&5	&7	&7	&7	&8	&9\\ \hline
			Type					&A	&B	&C	&B	&A	&B	&A	&B	&A	&A	&A	&B	&A	&A	&A	&A	&A\\ \hline
		\end{tabular}
		\end{table}
	\item Type $(2, 3, 5)\quad (p \neq 2, 3, 5)$ case:
		then $\frac{\beta_1}{\alpha_1} = \frac{1}{2}, \frac{\beta_2}{\alpha_2} \le \frac{2}{3}$ and $\frac{\alpha'}{\alpha_3} \ge \frac{6}{5}$.
		If $r = \lceil\frac{\beta_2}{\alpha_2} t + 1\rceil$,
			$s - r  - \frac{1}{2} t \ge \frac{6}{5} t + \frac{1}{5} - (\frac{2}{3} t + \frac{5}{3}) - \frac{1}{2} t = \frac{1}{30} t - \frac{22}{15}$.
		So (\ref{judge}) holds by a type A coboundary if $t \ge 44$.
		(\ref{judge}) also holds for $0 \le t \le 43$ cases by Table \ref{case3table}.
		\begin{table}[h]
		\begin{flushleft}
		\caption{Types of coboundaries in case (3).}
		\label{case3table}
		\begin{tabular}{@{}|c|c|c|c|c|c|c|c|c|c|c|c|c|c|c|c|c|c|}\hline
			$t$					&0	&1	&2	&3	&4	&5	&6	&7	&8	&9	&10	&11	&12	&13	&14\\ \hline
			$\mathrm{min}\{s\}$	&1	&2	&3	&4	&5	&7	&8	&9	&10	&11	&13	&14	&15	&16	&17\\ \hline	
			$r$					&1	&2	&2	&3	&0	&5	&5	&6	&6	&7	&8	&9	&9	&10	&10\\ \hline
			$\mathrm{min}\{s - r\}$	&0	&0	&1	&1	&5	&2	&3	&3	&4	&4	&5	&5	&6	&6	&7\\ \hline
			Type					&A	&B	&C	&B	&D	&B	&A	&B	&C	&B	&A	&B	&A	&B	&C\\ \hline
		\end{tabular}
		\begin{tabular}{@{}|c|c|c|c|c|c|c|c|c|c|c|c|c|c|c|c|c|c|}\hline
			$t$					&15	&16	&17	&18	&19	&20	&21	&22	&23	&24	&25	&26	&27	&28	&29\\ \hline
			$\mathrm{min}\{s\}$	&19	&20	&21	&22	&23	&25	&26	&27	&28	&29	&31	&32	&33	&34	&35\\ \hline	
			$r$					&11	&12	&13	&13	&14	&15	&15	&16	&17	&17	&18	&19	&19	&20	&21\\ \hline
			$\mathrm{min}\{s - r\}$	&8	&8	&8	&9	&9	&10	&11	&11	&11	&12	&13	&13	&14	&14	&14\\ \hline
			Type					&A	&A	&B	&A	&B	&A	&A	&A	&B	&A	&A	&A	&A	&A	&B\\ \hline
		\end{tabular}
		\begin{tabular}{@{}|c|c|c|c|c|c|c|c|c|c|c|c|c|c|c|c|c|}\hline
			$t$					&30	&31	&32	&33	&34	&35	&36	&37	&38	&39	&40	&41	&42	&43\\ \hline
			$\mathrm{min}\{s\}$	&37	&38	&39	&40	&41	&43	&44	&45	&46	&47	&49	&50	&51	&52\\ \hline	
			$r$					&21	&22	&23	&23	&24	&25	&25	&26	&27	&27	&28	&29	&29	&30\\ \hline
			$\mathrm{min}\{s - r\}$	&16	&16	&16	&17	&17	&18	&19	&19	&19	&20	&21	&21	&22	&22\\ \hline
			Type					&A	&A	&A	&A	&A	&A	&A	&A	&A	&A	&A	&A	&A	&A\\ \hline
		\end{tabular}
		\end{flushleft}
		\end{table}
\end{enumerate}

Then we have shown all monomial cocycles in $(x_0, y_0)$ are coboundaries.
All monomial terms in coordinate $(x_0', y_0')$ are sums of these terms and therefore coboundaries.
Subtracting these new coboundaries from coboundaries (\ref{polecoboundaryx}) and (\ref{polecoboundaryy}),
	cocycles (\ref{poleterm1}) and (\ref{poleterm2}) with a pole at $\{x_0 = 1\}$ are all coboundaries.
Then we have $H^1(\Theta_{P}) = 0$ and the proof has finished.

\subsection{Remarks on the proof.}
In the proof above, coboundaries of types other than A are all necessary for cohomology vanishing.
To get these coboundaries, all characteristics conditions are used.
In fact observing the list by Artin \cite{Art}, tautness and $F$-regularity are equivalent for rational double points.

On the other hand, there are some cases in each type of the star-shaped graphs whose cohomology calculation was omitted.
For example if the self-intersection number $- b_0$ of the central curve is sufficiently small, tautness holds for all characteristics
	because the $t = 0$ case calculation always holds and $\frac{\alpha'}{\alpha_3} \ge b_0 - 1$.
Even in the case $b_0 = 2$ and type $(2, 3, 5)$, which is the hardest case to vanish cohomology, if each branch has only one irreducible curve of self intersection $d_i$,
	then $\frac{\beta_1}{\alpha_1} = \frac{1}{2}, \frac{\beta_2}{\alpha_2} = \frac{1}{3}$ and $\frac{\alpha'}{\alpha_3} = \frac{9}{5}$
	and type B, C, D coboundaries are not necessary.
So $F$-regularity is not necessary for tautness in general.

\section{Discussions on $F$-pure rational cases}
\subsection{Relations between $F$-purity and tautness}

By Section 4.3, further relationships between $F$-singularity and tautness can be expected.
We discuss whether tautness holds for $F$-pure rational singularities or not.
The classification of $F$-pure rational singularities by Hara says that there are $F$-pure RDPs which are not $F$-regular.
This shows $F$-purity is not a sufficient condition for tautness of a rational singularity.

\begin{figure}[h]
	\centering
	\includegraphics[height=3cm]{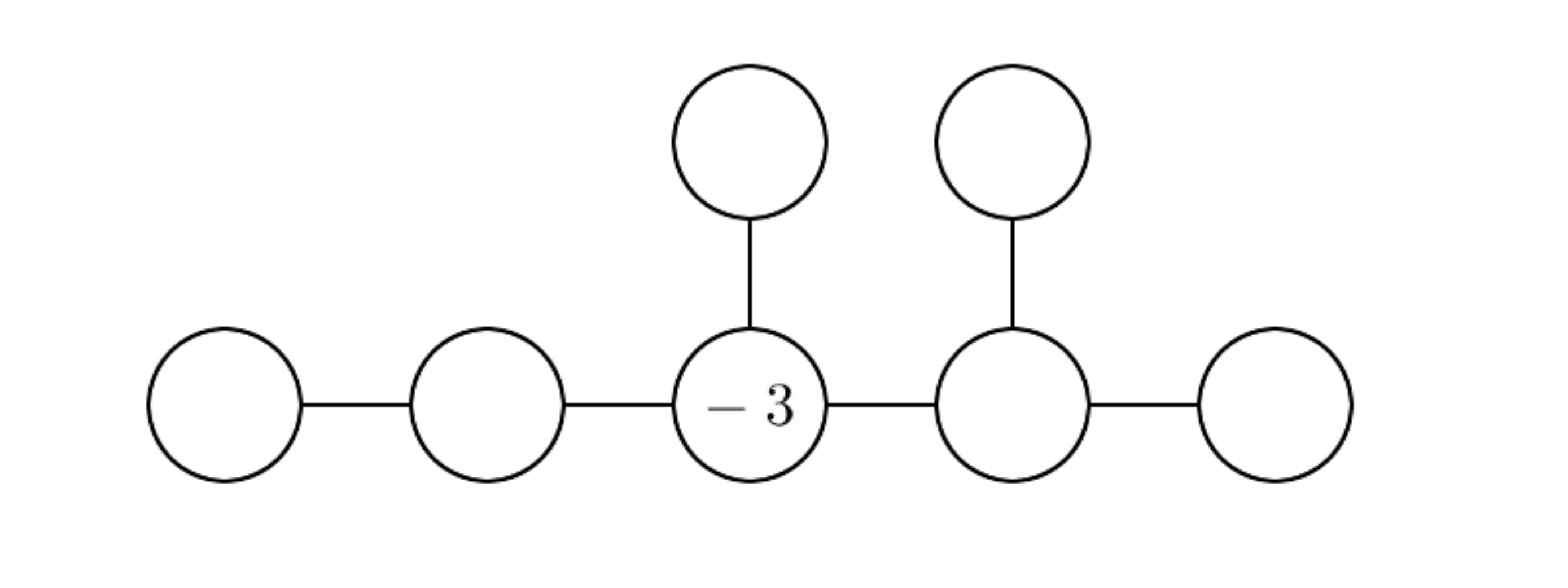}
	\caption{Non-$F$-pure rational taut graph.}
	\label{fig:counterex}
\end{figure}

On the other hand, the graph of a rational singularity shown in Figure \ref{fig:counterex} is a taut graph for large characteristics by \cite{Lau2} and \cite{Sch}
	but is not a graph of an $F$-pure singularity.
This means $F$-purity is not even a necessary condition for rational singularities to be taut.

Even though there is no implication between $F$-purity and tautness for rational singularities, some interesting phenomena can be observed.

\subsection{A kind of uniqueness for $F$-pure RDPs.}

In \cite{Art}, all rational double points in positive characteristics are presented using their defining equations as hypersurface singularities.
Using Fedder's criterion of $F$-purity \cite{Fed}, we can judge whether it is $F$-pure or not.
Results are shown in the following tables \ref{char2RDP}, \ref{char3RDP}, \ref{char5RDP} and \ref{char7RDP}.
Observing these tables, we can get the next theorem.

\begin{thm}
Let $(X, x)$ and $(X', x')$ be both two-dimensional $F$-pure rational double points over an algebraically closed field $k$ of a positive characteristic.
If $\Gamma_{X, x} \cong \Gamma_{X' x'}$, $(X, x)$ and $(X', x')$ are isomorphic to each other.
\end{thm}

\begin{table}[h]
	\caption{Rational double points in characteristic $2$.}
	\label{char2RDP}
	\begin{tabular}{|l|l|l|l|}\hline
		Graph				&Type							&Defining equation					&$F$-purity\\ \hline
		$A_n\, (n \ge 0)$		&$A_n$							&$z^{n + 1} + x y$					&$F$-pure\\ \hline
		$D_{2 n}\, (n \ge 2)$	&$D_{2 n}^0$						&$z^2 + x^2 y + x y^n$				&\\
							&$D_{2 n}^r\, (1 \ge r \ge n - 1)$		&$z^2 + x^2 y + x y^n + x y^{n - r} z$	&$\text{$F$-pure} \Leftrightarrow r = n - 1$\\ \hline
		$D_{2 n + 1}\, (n \ge 2)$	&$D_{2 n + 1}^0$					&$z^2 + x^2 y + y^n z$				&\\
							&$D_{2 n + 1}^r\, (1 \ge r \ge n - 1)$	&$z^2 + x^2 y + y^n z + x y^{n - r} z$	&$\text{$F$-pure} \Leftrightarrow r = n - 1$\\ \hline
		$E_6$				&$E_6^0$						&$z^2 + x^3 + y^2 z$				&\\
							&$E_6^1$						&$z^2 + x^3 + y^2 z + x y z$			&$F$-pure\\ \hline
		$E_7$				&$E_7^0$						&$z^2 + x^3 + x y^3$				&\\
							&$E_7^1$						&$z^2 + x^3 + x y^3 +x^2 y z$		&\\
							&$E_7^2$						&$z^2 + x^3 + x y^3 + y^3 z$			&\\
							&$E_7^3$						&$z^2 + x^3 + x y^3 + x y z$			&$F$-pure\\ \hline
		$E_8$				&$E_8^0$						&$z^2 + x^3 + y^5$					&\\
							&$E_8^1$						&$z^2 + x^3 + y^5 + x y^3 z$			&\\
							&$E_8^2$						&$z^2 + x^3 + y^5 + x y^2 z$			&\\
							&$E_8^3$						&$z^2 + x^3 + y^5 + y^3 z$			&\\
							&$E_8^4$						&$z^2 + x^3 + y^5 + x y z$			&$F$-pure\\ \hline
	\end{tabular}
\end{table}
\begin{table}[h]
	\caption{Rational double points in characteristic $3$.}
	\label{char3RDP}
	\begin{tabular}{|l|l|l|l|}\hline
		Graph			&Type		&Defining equation				&$F$-purity\\ \hline
		$A_n\, (n \ge 0)$	&$A_n$		&$z^{n + 1} + x y$				&$F$-pure\\ \hline
		$D_n\, (n \ge 4)$	&$D_n$		&$z^2 + x^2 y + y^{n - 1}$		&$F$-pure\\ \hline
		$E_6$			&$E_6^0$	&$z^2 + x^3 + y^4$				&\\
						&$E_6^1$	&$z^2 + x^3 + y^4 + x^2 y^2$	&$F$-pure\\ \hline
		$E_7$			&$E_7^0$	&$z^2 + x^3 + x y^3$			&\\
						&$E_7^1$	&$z^2 + x^3 + x y^3 + x^2 y^2$	&$F$-pure\\ \hline
		$E_8$			&$E_8^0$	&$z^2 + x^3 + y^5$				&\\
						&$E_8^1$	&$z^2 + x^3 + y^5 + x^2 y^3$	&\\
						&$E_8^2$	&$z^2 + x^3 + y^5 + x^2 y^2$	&$F$-pure\\ \hline
	\end{tabular}
\end{table}
\begin{table}[h]
	\caption{Rational double points in characteristic $5$.}
	\label{char5RDP}
	\begin{tabular}{|l|l|l|l|}\hline
		Graph			&Type		&Defining equation				&$F$-purity\\ \hline
		$A_n\, (n \ge 0)$	&$A_n$		&$z^{n + 1} + x y$				&$F$-pure\\ \hline
		$D_n\, (n \ge 4)$	&$D_n$		&$z^2 + x^2 y + y^{n - 1}$			&$F$-pure\\ \hline
		$E_6$			&$E_6$	&$z^2 + x^3 + y^4$					&$F$-pure\\ \hline
		$E_7$			&$E_7$	&$z^2 + x^3 + x y^3$				&$F$-pure\\ \hline
		$E_8$			&$E_8^0$	&$z^2 + x^3 + y^5$					&\\
						&$E_8^1$	&$z^2 + x^3 + y^5 + x y^4$			&$F$-pure\\ \hline
	\end{tabular}
\end{table}
\clearpage
\begin{table}[h]
	\caption{Rational double points in characteristic $\ge 7$.}
	\label{char7RDP}
	\begin{tabular}{|l|l|l|l|}\hline
		Graph			&Type		&Defining equation				&$F$-purity\\ \hline
		$A_n\, (n \ge 0)$	&$A_n$		&$z^{n + 1} + x y$				&$F$-pure\\ \hline
		$D_n\, (n \ge 4)$	&$D_n$		&$z^2 + x^2 y + y^{n - 1}$			&$F$-pure\\ \hline
		$E_6$			&$E_6$		&$z^2 + x^3 + y^4$				&$F$-pure\\ \hline
		$E_7$			&$E_7$		&$z^2 + x^3 + x y^3$			&$F$-pure\\ \hline
		$E_8$			&$E_8$		&$z^2 + x^3 + y^5$				&$F$-pure\\ \hline
	\end{tabular}
\end{table}

\subsection{Non-RDP star-shaped graphs with three branches.}
Next we see tautness of star-shaped graphs of non-RDP $F$-pure rational singularities with three branches, that is,
	the third cases of Hara's classification other than type $(2, 2, 2, 2)$.
We use the same method as Section 4.2.

First we show monomial form cocycles $x_0^s y_0^t y_0 \frac{\partial}{\partial y_0}$ are coboundaries.
For this, it was enough to show
\begin{equation}
	\lfloor\frac{\alpha'}{\alpha_3} t\rfloor + 1 - \lceil\frac{\beta_2}{\alpha_2} t\rceil - \frac{\beta_1}{\alpha_1} t \ge 0. \label{ytermjudge}
\end{equation}		
Next we show that monomial cocycles $x_0^s y_0^t x_0 \frac{\partial}{\partial x_0}$ are coboundaries using coboundaries of type A, B, C and D.
Recall that it is a coboundary from $U_1$ if $s \le \frac{\alpha'}{\alpha_3} t$.
Otherwise, it was enough to show there exists an integer $r$ that
	\begin{equation}
		(x_0 - 1)^r x_0^{s - r} y_0^t x_0 \frac{\partial}{\partial x_0} \label{xterm'}
	\end{equation}
	is a coboundary.

\begin{enumerate}
	\item Type $(3, 3, 3)\quad (p \equiv 1\, (\mathrm{mod}\, 3))$ case:
		since the intersection matrix is not negative definite if all self-intersections of irreducible curves are $- 2$ \cite{Bad},
			at least one component has self-intersection $- 3$ or less.
		So we may assume $\frac{\beta_1}{\alpha_1}, \frac{\beta_2}{\alpha_2} \le \frac{2}{3}$ and $\frac{\alpha'}{\alpha_3} \ge \frac{5}{3}$.
		
		If $t = 0$, then $\lfloor\frac{5}{3} t\rfloor + 1 - \lceil\frac{2}{3} t\rceil - \frac{2}{3} t = 1 \ge 0$.
		If $t \ge 1$, $\lfloor\frac{5}{3} t\rfloor + 1 - \lceil\frac{2}{3} t\rceil - \frac{2}{3} t
			\ge (\frac{5}{3} t + \frac{1}{3}) - (\frac{2}{3} t + \frac{2}{3}) - \frac{2}{3} t
			= \frac{1}{3} t - \frac{1}{3}
			\ge 0$.
		Therefore $x_0^s y_0^t y_0 \frac{\partial}{\partial y_0}$ terms are all coboundaries.
		
		Next we check cocycles $x_0^s y_0^t x_0 \frac{\partial}{\partial x_0}$.
		If $t \equiv 0\, (\mathrm{mod}\, 3)$, then (\ref{xterm'}) is a coboundary of type A if $r = \frac{2}{3} t + 1$.
		If $t \equiv 1\, (\mathrm{mod}\, 3)$, then (\ref{xterm'}) is a coboundary of type D if $r = 0$.
		If $t \equiv 2\, (\mathrm{mod}\, 3)$, then (\ref{xterm'}) is a coboundary of type C if $r = \frac{2}{3} t + \frac{2}{3}$.
		So $x_0^s y_0^t x_0 \frac{\partial}{\partial x_0}$ terms are all coboundaries.
		
	\item Type $(2, 3, 6)\quad (p \equiv 1\, (\mathrm{mod}\, 3))$ case:
		at least one self-intersection number is $- 3$ or less by the same reason as above.
		Then $\frac{\beta_1}{\alpha_1} = \frac{1}{2}, \frac{\beta_2}{\alpha_2} \le \frac{2}{3}, \frac{\alpha'}{\alpha_3} \ge \frac{7}{6}$
			and at least one inequality is not an equality.
		
		First consider the case $\frac{\beta_2}{\alpha_2} = \frac{1}{3}$.
		Since $\lfloor\frac{7}{6} t \rfloor + 1 - \lceil\frac{1}{3} t\rceil - \frac{1}{2} t \ge \frac{1}{3} t - \frac{1}{2}$,
			(\ref{ytermjudge}) holds for $t \ge 2$.
		Direct calculation shows that (\ref{ytermjudge}) also holds in $t = 0, 1$ cases.
		Next we consider the case $\frac{\alpha'}{\alpha_3} > \frac{7}{6}$, equivalently $\frac{\alpha'}{\alpha_3} \ge \frac{11}{6}$.
		Since $\lfloor\frac{11}{6} t\rfloor + 1 - \lceil\frac{2}{3} t\rceil - \frac{1}{2} t \ge \frac{2}{3} t - \frac{1}{2}$,
			(\ref{ytermjudge}) holds for $t \ge 1$.
		Direct calculation shows that (\ref{ytermjudge}) also holds in $t = 0$ case.
		Therefore cocycles $x_0^s y_0^t y_0 \frac{\partial}{\partial y_0}$ are all coboundaries.
		
		Next we check cocycles $x_0^s y_0^t x_0 \frac{\partial}{\partial x_0}$.
		If $\frac{\beta_2}{\alpha_2} = \frac{1}{3}$, then
			$\lfloor\frac{7}{6} t\rfloor + 1 - (\lceil\frac{1}{3} t\rceil + 1) - \frac{1}{2} \ge \frac{1}{3} t - \frac{3}{2}$
			implies $x_0^s y_0^t x_0 \frac{\partial}{\partial x_0}$ is a coboundary if $t \ge 5$.
		For $0 \le t \le 4$, Table \ref{236table1} gives the desired coboundaries.
		\begin{center}
		\label{236table1}
		\begin{tabular}{@{}|c|c|c|c|c|c|c|c|c|c|c|c|c|c|c|c|c|}\hline
			$t$					&0	&1	&2	&3	&4\\ \hline
			$\mathrm{min}\{s\}$		&1	&2	&3	&4	&5\\ \hline	
			$r$					&1	&2	&2	&2	&3\\ \hline
			$\mathrm{min}\{s - r\}$	&0	&0	&1	&2	&2\\ \hline
			Type					&A	&B	&A	&A	&A\\ \hline
		\end{tabular}
		\end{center}
		If $\frac{\alpha'}{\alpha_3} = \frac{11}{6}$, then
			$\lfloor\frac{11}{6} t\rfloor + 1 - (\lceil\frac{2}{3} t\rceil + 1) - \frac{1}{2} t \ge \frac{2}{3} t - \frac{3}{2}$
			says that $x_0^s y_0^t x_0 \frac{\partial}{\partial x_0}$ is a coboundary if $t \ge 3$.
		If $t = 0$, $(x_0 - 1) x_0 \frac{\partial}{\partial x_0}$ is a type A coboundary.
		If $t = 1$, $(x_0 - 1)^2 y_0 x_0 \frac{\partial}{\partial x_0}$ is a type B coboundary.
		If $t = 2$, $(x_0 - 1)^2 x_0 y_0^2 x_0 \frac{\partial}{\partial x_0}$ is a type C coboundary.
		Therefore cocycles $x_0^s y_0^t x_0 \frac{\partial}{\partial x_0}$ are all coboundaries.
		
	\item Type $(2, 4, 4)\quad (p \equiv 1\, (\mathrm{mod}\, 4))$ case:
		same argument as above shows at least one component has self-intersection $- 3$ or less.
		So we may assume $\frac{\beta_1}{\alpha_1} = \frac{1}{2},\, \frac{\beta_2}{\alpha_2} = \frac{3}{4}$ and $\frac{\alpha'}{\alpha_3} \ge \frac{7}{4}$.
		
		Since $\lfloor\frac{7}{4} t \rfloor + 1 - \lceil\frac{3}{4} t\rceil - \frac{1}{2} t \ge \frac{1}{2} t - \frac{1}{2}$,
			(\ref{ytermjudge}) holds for $t \ge 1$.
		$t = 0$ case is calculated directly and (\ref{ytermjudge}) holds.
		
		Next check $x_0^s y_0^t x_0 \frac{\partial}{\partial x_0}$ is a coboundary.
		It is a coboundary if $t \ge 3$ because $\lfloor\frac{7}{4} t\rfloor + 1 - (\lfloor\frac{3}{4} t\rfloor + 1) - \frac{1}{2} t \ge \frac{1}{2} t - \frac{3}{2}$.
		If $t = 0$, $(x_0 - 1) x_0 \frac{\partial}{\partial x_0}$ is a type A coboundary.
		If $t = 1$, $(x_0 - 1)^2 y_0 x_0 \frac{\partial}{\partial x_0}$ is a type B coboundary.
		If $t = 2$, $(x_0 - 1)^3 x_0 y_0^2 x_0 \frac{\partial}{\partial x_0}$ is a type A coboundary.
		This shows that cocycles $x_0^s y_0^t x_0 \frac{\partial}{\partial x_0}$ are all coboundaries.
\end{enumerate}

As a result, $H^1(\Theta_P) = 0$ is shown in these cases.

\subsection{Type $(2, 2, 2, 2)$ star-shaped graphs.}
Though it was seen in Section 5.2 that $F$-purity does not implies tautness of rational singularities, there still remains a possibility that Theorem 5.1 holds for non-RDPs.
Unfortunately, star-shaped graphs of type $(2, 2, 2, 2)$ is a counterexample of this, not only of tautness.

Fix a graph $\Gamma$ of type $(2, 2, 2, 2)$ appearing as a graph of an $F$-pure rational singularity.
It was shown that there are infinitely many $\lambda$s which gives the intersection points at the central curve $0, - 1, \lambda, \infty \in \mathbb{P}_k^1$.
The permutation group $\mathfrak{S}_4$ acts on these $\lambda$ and different orbits represent different positions of intersections.
Since there are infinitely many orbits, there are infinite family of exceptional curves $\{E_{\bar{\lambda}}\}_{\bar{\lambda}}$
	embedded in nonsingular surfaces and associated with $\Gamma$.
These $E_{\bar{\lambda}}$ always satisfy the condition of contractibility \cite{Bad}, they can be contracted into rational singularities.
Then we obtain infinitely many non-isomorphic $F$-pure rational singularities whose graphs are all $\Gamma$.
This gives a counterexample of Theorem 5.1 in non-RDP case.

\subsection{$_*\tilde{D}_{n + 3} (n \ge 2)$ graphs}
In the case $k = \mathbb{C}$, a $_*\tilde{D}_{n + 3} (n \ge 2)$ graph is always taut if the negativity of the intersection matrix is satisfied \cite{Lau1}.
In arbitrary positive characteristic, there are examples of $_*\tilde{D}_{n + 3} (n \ge 2)$ graphs which might not be taut.
This is because $H^1(\Theta_P) \neq 0$, but we need to improve Sch\"{u}ller's criterion to judge whether it is taut or not in fact.

\begin{figure}[h]
	\centering
	\includegraphics[height=3cm]{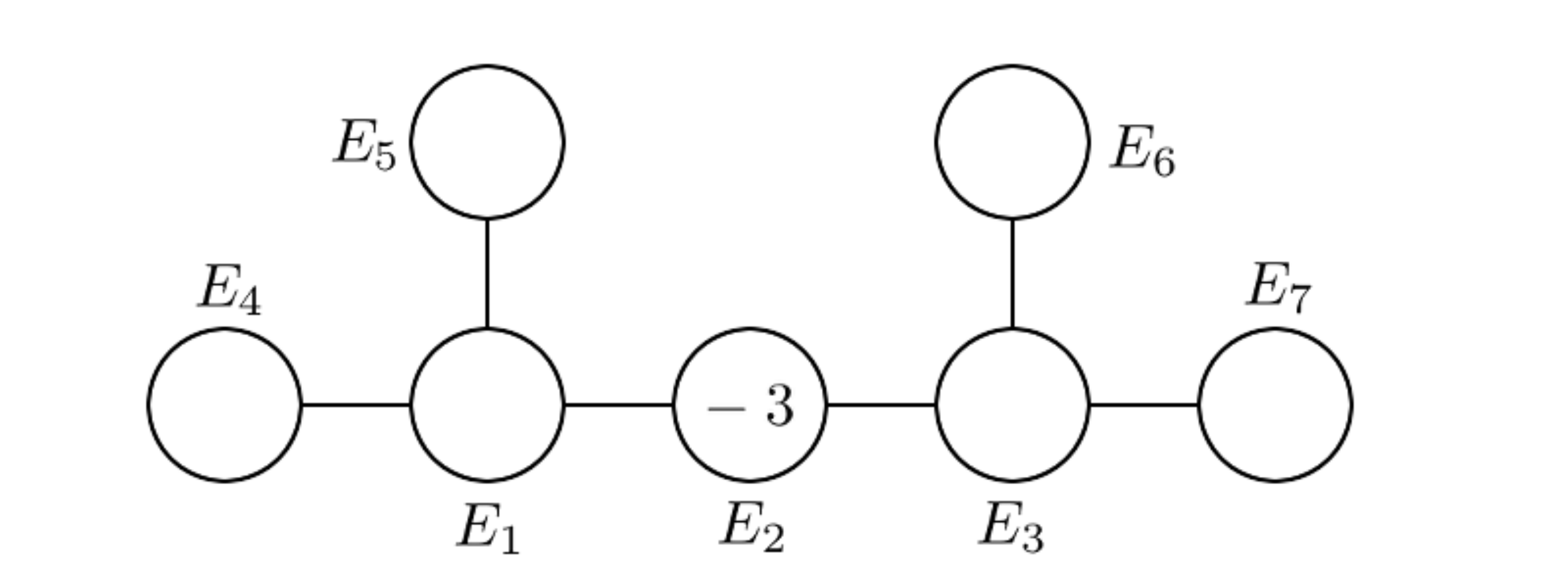}
	\caption{An example of $_*\tilde{D}_{n + 3} (n \ge 2)$ graph with $H^1(\Theta_P) \neq 0$ in $p = 3$.}
	\label{fig:Dn+3ex}
\end{figure}

We see one example.
If $p = 3$, the graph in Figure \ref{fig:Dn+3ex} gives $P$ such that $H^1(\Theta_P) \neq 0$.
Then we can take $\tilde{\nu} = 4,\, \nu_1 = \nu_3 = 28,\, \nu_2 = 20$ and $\nu_4 = \nu_5 = \nu_6 = \nu_7 = 16$.

We can calculate $H^1(\Theta_P)$ using a covering similar to the one in Section 4.
That is, we see four vertices in the left as a subgraph of a star-shaped graph of type $(2, 2, d)$ and define $U_0$ in the same way.
Then $U_0 = (E_1 \setminus E_2) \cup E_4 \cup E_5$.
Set $U_2 = (E_3 \setminus E_2) \cup E_6 \cup E_7$ for the opposite side in the same way.

Then $\Gamma(\Theta_{U_0})$ has a basis as follows:
\begin{align}
	&(x_0 - 1)^r x_0^s y_0^t x_0 \frac{\partial}{\partial x_0}\\
		&\left(t \ge 0,\, r = \mathrm{min}\left\{\lceil\frac{1}{2} t\rceil + 1, 16\right\},\,
			s \ge \mathrm{min}\left\{\frac{1}{2} t, 15\right\}\right)\notag\\
	&(x_0 - 1)^r x_0^s y_0^t y_0 \frac{\partial}{\partial y_0}\\
		&\left(t \ge 0,\, r = \mathrm{min}\left\{\lceil\frac{1}{2} t\rceil, 16\right\},\,
			s \ge \mathrm{min}\left\{\frac{1}{2} t, 16\right\}\right)\notag\\
	&(x_0 - 1)^r x_0^s y_0^t (2 x_0 \frac{\partial}{\partial x_0} - y_0 \frac{\partial}{\partial y_0})\\
			&\left(r, t \ge 1,\, r = \mathrm{min}\left\{\frac{1}{2} t + \frac{3}{2}, 16\right\},\, s = \frac{1}{2} t - \frac{1}{2}\right)\notag\\
	&2 (x_0 - 1)^{r + 1} x_0^s y_0^t x_0 \frac{\partial}{\partial x_0} - (x_0 - 1)^r x_0^s y_0 \frac{\partial}{\partial y_0}\\
			&\left(r, t \ge 1,\, r = \frac{1}{2} t - \frac{1}{2},\, s = \mathrm{min}\left\{\frac{1}{2} t + \frac{1}{2}, 16\right\}\right).\notag
\end{align}

Let $U_1$ be $(E_1 \setminus (E_4 \cup E_5)) \cup E_2 \cup (E_3 \setminus (E_6 \cup E_7))$.
Then $P = U_0 \cup U_1 \cup U_2$ is a Leray cover for $\Theta_P$ provided $H^1(\Theta_{U_1}) = 0$.

First compute $H^1(\Theta_{U_1}) = 0$ and a basis of $\Gamma(\Theta_{U_1})$.
Let $A_1 = (E_1 \setminus (E_4 \cup E_5)) \cup (E_2 \setminus E_3)$ and $A_2 = (E_2 \setminus E_1) \cup (E_3 \setminus (E_6 \cup E_7))$.
Then $U_1 = A_1 \cup A_2$ is an affine covering.
We set the coordinates as
\begin{align*}
	&A_1 = \mathrm{Spec}(k[x_1, y_1, (1 - x_1)^{- 1}] / (x_1^{20} y_1^{28})\\
	&A_2 = \mathrm{Spec}(k[x_2, y_2, (1 - y_2)^{- 1}] / (x_2^{28} y_2^{20})\\
	&x_1 = x_2^3 y_2,\, y_1 = x_2^{- 1}.
\end{align*}
Then
\begin{equation*}
	x_1 \frac{\partial}{\partial x_1} = y_2 \frac{\partial}{\partial y_2},\,
		y_1 \frac{\partial}{\partial y_1} = - x_2 \frac{\partial}{\partial x_2} + 3 y_2 \frac{\partial}{\partial y_2} = - x_2 \frac{\partial}{\partial x_2}.
\end{equation*}
Here vanishing of the $y_2 \frac{\partial}{\partial y_2}$ term in the second equation by the characteristic condition is the key point.
By this formula,
	\begin{align*}
		x_1^s y_1^t x_1 \frac{\partial}{\partial x_1} &= x_2^{3 s - t} y_2^s y_2 \frac{\partial}{\partial y_2},\\
		x_1^s y_1^t y_1 \frac{\partial}{\partial y_1} &= - x_2^{3 s - t} y_2^s x_2 \frac{\partial}{\partial x_2}.
	\end{align*}
This situation is similar to the proof of lemma 4.2 and $H^1(\Theta_{U_1}) = 0$ can be checked in the same way.

By the coordinate change given above,
	$x_1^s y_1^t x_1 \frac{\partial}{\partial x_1} \in \Gamma(\Theta_{U_1})
		\Leftrightarrow x_1^s y_1^t y_1 \frac{\partial}{\partial y_1} \in \Gamma(\Theta_{U_1})
		\Leftrightarrow 3 s \ge t$
	for small $s$.
Precisely, $\Gamma(\Theta_{U_1})$ has a basis as follows:
\begin{align*}
	x_1^s y_1^t x_1 \frac{\partial}{\partial x_1} = x_2^{3 s - t} y_2^s y_2 \frac{\partial}{\partial y_2}&\quad (t \ge 0, \frac{1}{3} t \le s \le 18)\\
	x_1^s y_1^t y_1 \frac{\partial}{\partial y_1} = - x_2^{3 s - t} y_2^s x_2 \frac{\partial}{\partial x_2}&\quad (t \ge 0, \frac{1}{3} t \le s \le 19)\\
	x_1^s y_1^t x_1 \frac{\partial}{\partial x_1}&\quad (t \ge 0, s \ge 19)\\
	x_1^s y_1^t y_1 \frac{\partial}{\partial y_1}&\quad (t \ge 0, s \ge 20)\\
	x_2^s y_2^t x_2 \frac{\partial}{\partial x_2}&\quad (s \ge 0, t \ge 20)\\
	x_2^s y_2^t y_2 \frac{\partial}{\partial y_2}&\quad (s \ge 0, t \ge 19)\\
	\left((1 - x_1)^r - (1 + x_1 + \cdots + x_1^{18})^{- r}\right) &y_1^t x_1 \frac{\partial}{\partial x_1}\quad (r < 0,\, t \ge 0)\\
	\left((1 - x_1)^r - (1 + x_1 + \cdots + x_1^{19})^{- r}\right) &y_1^t y_1 \frac{\partial}{\partial y_1}\quad (r < 0,\, t \ge 0)\\
	\left((1 - y_2)^r - (1 + y_2 + \cdots + y_2^{19})^{- r}\right) &x_2^s x_2 \frac{\partial}{\partial x_2}\quad (r < 0,\, t \ge 0)\\
	\left((1 - y_2)^r - (1 + y_2 + \cdots + y_2^{18})^{- r}\right) &x_2^s y_2 \frac{\partial}{\partial y_2}\quad (r < 0,\, t \ge 0).
\end{align*}

Since what we want is only their images in $\Gamma(\Theta_{U_0 \cap U_1})$ and $\Gamma(\Theta_{U_1 \cap U_2})$,
	we use the expression by $(x_0, y_0)$ and $(x_3, y_3)$.
Note that $\Gamma(\Theta_{U_0 \cap U_1})$ has a basis as follows:
\begin{equation}
	\begin{aligned}
		x_0^s y_0^t x_0 \frac{\partial}{\partial x_0}&\quad (s \in \mathbb{Z}, 0 \le t \le 27)\\
		(x_0 - 1)^r x_0^{- r} x_0^{2 t} y_0^t x_0 \frac{\partial}{\partial x_0}&\quad (r < 0, 0 \le t \le 27)\\
		x_0^s y_0^t y_0 \frac{\partial}{\partial y_0}&\quad (s \in \mathbb{Z}, 0 \le t \le 26)\\
		(x_0 - 1)^r x_0^{- r} x_0^{2 t} y_0^t y_0 \frac{\partial}{\partial y_0}&\quad (r < 0, 0 \le t \le 26)
	\end{aligned} \label{basiscoc0}
\end{equation}
This is different from the one used in Section 4, but convenient in this case.

Via the coordinate change,
	$\mathrm{Im}\left(\Gamma(\Theta_{U_1}) \rightarrow \Gamma(\Theta_{U_0 \cap U_1}) \oplus \Gamma(\Theta_{U_1 \cap U_2})\right)$ has a basis as follows:
\begin{align}
	(x_0^s y_0^t x_0 \frac{\partial}{\partial x_0},\, x_3^{8 t - 5 s} y_3^{5 t - 3 s} \left(x_3 \frac{\partial}{\partial x_3} - y_3 \frac{\partial}{\partial y_3}\right))
		\quad (0 \le &t \le 27, \frac{5}{3} t - 9 \le s \le \frac{5}{3} t) \label{cob:ordx}\\
	(x_0^s y_0^t y_0 \frac{\partial}{\partial y_0},\, - x_3^{8 t - 5 s} y_3^{5 t - 3 s}  y_3 \frac{\partial}{\partial y_3})
		\quad (0 \le &t \le 26, \frac{5}{3} t - \frac{26}{3} \le s \le \frac{5}{3} t)\label{cob:ordy}\\
	(x_0^s y_0^t x_0 \frac{\partial}{\partial x_0},\, 0)\quad (0 \le &t \le 27, s < \frac{5}{3} t - 9)\\
	(x_0^s y_0^t y_0 \frac{\partial}{\partial y_0},\, 0)\quad (0 \le &t \le 26, s < \frac{5}{3} t - \frac{26}{3})\\
	(0,\, x_3^s y_3^t x_3 \frac{\partial}{\partial x_3})\quad (0 \le &t \le 27, s < \frac{5}{3} t - 9)\\
	(0,\, x_3^s y_3^t y_3 \frac{\partial}{\partial y_3})\quad (0 \le &t \le 26, s < \frac{5}{3} t - \frac{26}{3})\\
	(\left(x_0^{- r} (x_0 - 1)^r - (1 + x_0^{- 1} + \cdots + x_0^{- 18})^{- r}\right) &x_0^{2 t} y_0^t x_0 \frac{\partial}{\partial x_0},\, 0)\\
		&\quad (r < 0,\, 0 \le t \le 27)\notag\\
	(\left(x_0^{- r} (x_0 - 1)^r - (1 + x_0^{- 1} + \cdots + x_0^{- 18})^{- r}\right) &x_0^{2 t} y_0^t y_0 \frac{\partial}{\partial y_0},\, 0) \label{cob:vanpoley0}\\
		&\quad (r < 0,\, 0 \le t \le 26)\notag\\
	(0,\, \left(x_3^{- r} (x_3 - 1)^r - (1 + x_3^{- 1} + \cdots + x_3^{- 18})^{- r}\right) &x_3^{2 t} y_3^t x_3 \frac{\partial}{\partial x_3}) \label{cob:vanpolex3}\\
		&\quad (r < 0,\, 0 \le t \le 27)\notag\\
	(0,\, \left(x_3^{- r} (x_3 - 1)^r - (1 + x_3^{- 1} + \cdots + x_3^{- 18})^{- r}\right) &x_3^{2 t} y_3^t y_3 \frac{\partial}{\partial y_3}) \label{cob:vanpoley3}\\
		&\quad (r < 0,\, 0 \le t \le 26)\notag
\end{align}

We show that $(0, x_3 \frac{\partial}{\partial x_3}) \in \Gamma(\Theta_{U_0 \cap U_1}) \oplus \Gamma(\Theta_{U_1 \cap U_2})$ is not a coboundary.
First consider coboundaries from $U_2$.
In the basis of $\Gamma(\Theta_{U_2})$ above, coboundaries related to $x_3 \frac{\partial}{\partial x_3}$ are only
\begin{equation}
	(x_3 - 1) x_3^s x_3 \frac{\partial}{\partial x_3}\quad (s \ge 0). \label{cob:target}
\end{equation}
This gives relations between $x_3 \frac{\partial}{\partial x_3}$ and $x_3^s x_3 \frac{\partial}{\partial x_3}\quad (s \ge 1)$ and no relations to others.

We check coboundaries from $U_1$ related to elements in (\ref{cob:target}).
Observing (\ref{cob:ordx}), related elements satisfy $8 t - 5 s \ge 0$ and $5 t - 3 s = 0$.
This implies $t =0$ and $x_0 \frac{\partial}{\partial x_0} = x_3 \frac{\partial}{\partial x_3} - y_3 \frac{\partial}{\partial y_3}$
	is the only coboundary we can use to vanish the target cocycle.
Next see coboundaries of the form (\ref{cob:vanpolex3}).
Related terms satisfy $t = 0$.
On the other hand, $x_3^{- r} (x_3 - 1)^r x_3 \frac{\partial}{\partial x_3}$ has a pole of order $- r$ at $x_3 = 1$
	and therefore no coboundaries from $U_2$ have terms to cancel this pole.
This tells us that any coboundaries including a nontrivial sum of these with $r < 0$ always have a pole at $x_3 = 1$ and we cannot use them to make
	$(0, x_3 \frac{\partial}{\partial x_3})$.
Terms from (\ref{cob:ordy}) to (\ref{cob:vanpoley0}) and (\ref{cob:vanpoley3}) has no related terms.

We can apply the same argument done for $x_3 \frac{\partial}{\partial x_3}$ to $x_0 \frac{\partial}{\partial x_0}$ by symmetricity.
Consider the image of
\[(\theta_0, \xi x_0 \frac{\partial}{\partial x_0} + \theta_1, \theta_2) \in \Gamma(\Theta_{U_0}) \oplus \Gamma(\Theta_{U_1}) \oplus \Gamma(\Theta_{U_2})\quad (\xi \in k)\]
where $\theta_1$ does not have $x_0 \frac{\partial}{\partial x_0}$ term.
Then its image by the restriction map is
\begin{equation*}
	(\xi x_0 \frac{\partial}{\partial x_0} + \theta_0 + \theta_1, \xi x_0 \frac{\partial}{\partial x_0} - \xi y_0 \frac{\partial}{\partial y_0} + \theta_0 + \theta_2).
\end{equation*}
Here the signature multiplied to the restriction is set all positive because this change is not intrinsic.
Using the basis (\ref{basiscoc0}),
\[\xi x_0 \frac{\partial}{\partial x_0} + \theta_0 + \theta_1
	= \left(\sum_{s \in \mathbb{Z}} \xi_s x_0^s x_0 \frac{\partial}{\partial x_0}\right) +\, \text{other terms}\]
where $\sum_{s \ge 0} \xi_s = \xi$.
Similarly,
\[\theta_2 - \xi x_0 \frac{\partial}{\partial x_0} + \theta_0
	= \left(\sum_{s \in \mathbb{Z}} \xi_s' x_3^s x_3 \frac{\partial}{\partial x_3}\right) +\, \text{other terms}\]
where $\sum_{s \ge 0} \xi_s = \xi$.

$(0, x_3 \frac{\partial}{\partial x_3})$ satisfies $\sum_{s \ge 0} \xi_s = 0$ and $\sum_{s \ge 0} \xi_s = 1$.
This implies that it is not a coboundary and $H^1(\Theta_P) \neq 0$ is proved.

\vspace{1\baselineskip}
\textbf{Acknowledgements.}
This is a master course thesis in Graduate School of Mathematical Science, the University of Tokyo.
The author thanks all who helped his work.
Especially his advisor Shunsuke Takagi always gave me inspective advice in the weekly seminars and I could find the direction of my research.
Discussions with my colleagues at the University of Tokyo, especially frequent chats with Ippei Nagamachi and Hironori Matsuue helped him
	reorganize knowledge and find another point of view.
My family supported me in daily life and I could spend much time for my work.


\end{document}